\input ./style/arxiv-general.cfg
\documentclass[aap,MSNbibl,seceqn,dvips]{arximspdf}
\makeatletter
   \@ifpackageloaded{graphicx}{}{\usepackage{graphicx}}
\makeatother

\doi{10.1214/14-AAP1094}
\volume{26}
\issue{1}
\pubyear{2016}
\firstpage{360}
\lastpage{401}
\docsubty{FLA}

\makeatletter
\newcommand{\rd}{\mathrm{d}}
\newcommand{\rrvert}{\vert}
\newcommand{\rrVert}{\Vert}
\newcommand{\llvert}{\vert}
\newcommand{\llVert}{\Vert}
\renewcommand{\mid}{|}
\newtheorem{teo}{Theorem}[section]
\newtheorem{cor}{Corollary}[section]
\newtheorem{prop}{Proposition}[section]
\newtheorem{lemma}{Lemma}[section]
\newproclaim{remark}{Remark}[section]
\newcommand{\ve}{\varepsilon}
\newcommand{\gep}{\epsilon}
\newcommand{\rmd}{\mathrm{d}}
\newcommand{\rmi}{\mathrm{i}}
\newcommand{\whH}{\widehat{H}}
\newcommand{\whV}{\widehat{V}}
\newcommand{\whL}{\widehat{L}}
\newcommand{\whX}{\widehat{X}}
\newcommand{\whk}{\widehat{\kappa}}
\newcommand{\Fbar}{\overline{F}}
\newcommand{\pibar}{\overline{\Pi}}
\newcommand{\Vbar}{\overline{V}}
\newcommand{\Xbar}{\overline{X}}
\newcommand{\ceb}{\overline{\mathcal E}}
\newcommand{\nbar}{\overline{n}}
\newcommand{\Ebar}{\overline{E}}
\newcommand{\R}{\mathbf{R}}
\newcommand{\ga}{\alpha}
\newcommand{\gl}{\lambda}
\newcommand{\gk}{\kappa}
\newcommand{\gz}{\zeta}
\newcommand{\gb}{\beta}
\newcommand{\gt}{\theta}
\newcommand{\gd}{\delta}
\newcommand{\gs}{\sigma}
\newcommand{\tf}{\tilde f}
\newcommand{\tn}{\tilde n}
\newcommand{\tj}{\tilde j}
\newcommand{\rme}{\mathbf{e}}
\newcommand{\vep}{\varepsilon}

\makeatother

\begin{document}
\begin{frontmatter}

\title{Sample path behavior of a L\'evy insurance risk process
approaching ruin, under the Cram\'er--Lundberg and convolution
equivalent~conditions\thanksref{T1}}
\runtitle{L\'evy insurance risk process}

\begin{aug}
\author[A]{\fnms{Philip S.}~\snm{Griffin}\corref{}\ead[label=e1]{psgriffi@syr.edu}}
\runauthor{P.~S. Griffin}
\affiliation{Syracuse University}
\address[A]{Department of Mathematics\\
Syracuse University\\
Syracuse, New York 13244-1150\\
USA\\
\printead{e1}}
\end{aug}
\thankstext{T1}{Supported in part by Simons Foundation Grant 226863.}

%
\received{\smonth{9} \syear{2013}}
%
\revised{\smonth{10} \syear{2014}}

%
\begin{abstract}
Recent studies have demonstrated an interesting connection between
the asymptotic behavior at ruin of a L\'evy insurance risk process
under the Cram\'er--Lundberg
and convolution equivalent conditions.
For example, the limiting distributions of the overshoot and the
undershoot are strikingly similar in these two settings. This is
somewhat surprising since the global sample path behavior of the
process under these two conditions is quite different. Using tools
from excursion theory and fluctuation theory, we provide a means
of transferring results from one setting to the other which, among
other things, explains this connection and leads to new asymptotic
results. This is done by describing the evolution
of the sample paths from the time of the last maximum prior to ruin
until ruin occurs.
\end{abstract}

%
\end{frontmatter}

\setcounter{footnote}{1}
\section{Introduction}\label{s1}\label{sec1}

It is becoming increasingly popular to model insurance risk processes
with a general L\'evy process.
In addition to new and interesting mathematics, this approach allows
for direct modeling of aggregate claims
which can then be calibrated against real aggregate data, as opposed to
the traditional approach of modeling individual claims. Whether this
approach is superior remains to be seen, but it offers, at a minimum,
an alternative, to the traditional approach.
The focus of this paper will be on two such L\'evy models,
and their sample path behavior as ruin approaches.

Let $X=\{X_{t}\dvtx  t \geq0 \}$,
$X_0=0$, be a L\'{e}vy process
with characteristics $(\gamma, \sigma^2, \Pi_X)$.
The characteristic function of $X$
is given by the L\'{e}vy--Khintchine
representation, $Ee^{\mathrm{i}\theta X_{t}} = e^{t \Psi_X(\theta)}$,
where
\[
\label{lrep} \Psi_X(\theta) = \rmi\theta\gamma-
\sigma^2\theta^2/2+ \int_{\R}
\bigl(e^{\rmi\theta x}-1- \rmi\theta x \mathbf{1}_{\{\llvert  x\rrvert  <1\}}\bigr)
\Pi_X(d x)\qquad\mbox{for } \theta\in\mathbf{R}.
\]
To avoid trivialities, we assume $X$ is nonconstant.
In the insurance risk model, $X$ represents the excess in
claims over premium. An insurance company starts with
an initial positive reserve
$u$, and ruin occurs
if this level is exceeded by $X$. To reflect the insurance company's
desire to collect sufficient premia to prevent almost certain ruin, it
is assumed that $X_t\to-\infty$ a.s.
This is the \textit{general L\'evy insurance risk model}, which we will
investigate under two distinct conditions. The first is the well-known
Cram\'er--Lundberg condition:
%
\begin{equation}
\label{Cr} Ee^{\ga X_1}=1\quad\mbox{and}\quad EX_1e^{\ga X_1}<
\infty\qquad \mbox {for some } \ga>0.
\end{equation}
The second, introduced by Kl\"{u}ppelberg, Kyprianou and Maller \cite
{kkm}, is the convolution equivalent condition:
%
\begin{equation}
\label{Sa} Ee^{\ga X_1}<1\quad\mbox{and}\quad X_1^+\in{
\mathcal S}^{(\alpha)}  \qquad\mbox{for some } \ga>0,
\end{equation}
where ${\mathcal S}^{(\alpha)}$ denotes the class of convolution
equivalent distributions of index $\ga$. The formal description of
${\mathcal S}^{(\alpha)}$ will be given in Section~\ref{s6}. Typical
examples of distributions in ${\mathcal S}^{(\alpha)}$ are those with
tails of the form
\[
\label{p} P(X_1>x)\sim\frac{e^{-\alpha x}}{x^p} \qquad\mbox{for } p>1.
\]
Under (\ref{Sa}), $Ee^{\gt X_1}=\infty$ for all $\gt>\ga$, so (\ref
{Cr}) must fail. Hence, conditions (\ref{Cr}) and (\ref{Sa}) are
mutually exclusive. For a further comparison, see the introduction to~\cite{GM2}.

Historically, the first insurance risk model to be extensively studied
was the \textit{compound Poisson model}.\footnote{This is often called the
Cram\'er--Lundberg model, as opposed to the Cram\'er--Lundberg
condition~(\ref{Cr}).}
This arises when $X$ is a spectrally positive compound Poisson process
with negative drift.
In recent years, attention has turned to the general L\'evy insurance
risk model (see Kyprianou \cite{kypbook} for a detailed discussion of
the general model), where
considerable progress has been made in calculating the limiting
distribution of several variables related to ruin; see Doney, Kl\"
{u}ppelberg and Maller \cite{DKM}, Doney and Kyprianou \cite{DK},
Griffin and Maller \cite{GM2}, Kl\"{u}ppelberg, Kyprianou and Maller
\cite{kkm} and the references therein. To give some examples,
particularly relevant to this paper, we first need a little notation. Set
\begin{eqnarray*}
\label{Xbar} \Xbar_t &=&\sup_{0\le s\le t} X_s,
\\
\label{tu} \tau(u) &=&\inf\bigl\{t\dvtx X(t)>u\bigr\},
\end{eqnarray*}
and let $P^{(u)}$ denote the probability measure $P^{(u)}( \cdot )=P( \cdot \mid \tau(u)<\infty)$. Let
$H$ be the ascending ladder height process, and $\Pi_{H}, \rmd_H$ and
$q$ its L\'evy measure, drift and killing rate, respectively; see
Section~\ref{s2} for more details. Then under the Cram\'er--Lundberg
condition (\ref{Cr}), it was shown in \cite{GMS} that the limiting
distributions of the shortfall at, and the minimum surplus prior to,
ruin are given by
%
\begin{eqnarray}
\label{oCr} P^{(u)}(X_{\tau(u)}-u\in d x) & \stackrel{\mathrm
{w}} {\longrightarrow} &q^{-1}\ga \biggl[\rmd_H
\delta_{0}(d x)+ \int_{y\ge0} e^{\alpha y}
\Pi_{
H}(y+d x)\,d y \biggr],\hspace*{-25pt}
\nonumber\\[-8pt]\\[-8pt]\nonumber
P^{(u)}(u-\Xbar_ {{\tau(u)}-}\in d y)&\stackrel{\mathrm {w}} {
\longrightarrow} &q^{-1} \alpha\bigl[\rmd_H
\delta_{0}(d y)+ e^{\alpha y}\pibar_{H}(y)\,d y\bigr],
\end{eqnarray}
where $ \stackrel{\mathrm{w}}{\longrightarrow}$ denotes weak
convergence and $\delta_{0}$ is a point mass
at $0$.
Under the convolution equivalent condition (\ref{Sa}), it follows from
Theorem 4.2 in \cite{kkm} and Theorem~10 in \cite{DK} (see also
Section~7 of \cite{GM2}) that the corresponding limits are
%
\begin{eqnarray}\label{oSa}
P^{(u)}(X_{\tau(u)}-u\in d x)
&\stackrel{\mathrm{w}} {\longrightarrow}& q^{-1}\ga \biggl[-\ln\bigl(Ee^{\ga H_1}
\bigr)e^{-\alpha x}\,d x +\rmd_H\delta_{0}(d x)\nonumber
\\
&&\hspace*{58pt}{} + \int
_{y\ge0} e^{\alpha y} \Pi_{
H}(y+d x)\,d y \biggr],
\\
P^{(u)}(u-\Xbar_ {{\tau(u)}-}\in d y) &\stackrel{\mathrm {v}} {
\longrightarrow}&  q^{-1} \alpha\bigl[\rmd_H
\delta_{0}(d y) + e^{\alpha y}\pibar_{H}(y)\,d y\bigr],\nonumber
\end{eqnarray}
where $\stackrel{\mathrm{v}}{\longrightarrow}$ denotes vague
convergence of measures on $[0,\infty
)$.\footnote{In (\ref{oCr}) and (\ref{oSa}), it is assumed that $X_1$
has a nonlattice distribution. Similar results hold in the lattice case
if the limit is taken through points in the lattice span, but to avoid
repetition we will henceforth make the nonlattice assumption.}

The resemblance between the results in (\ref{oCr}) and (\ref{oSa}) is
striking and in many ways quite surprising, since the paths resulting
in ruin behave very differently in the two cases as we now explain.
Under the Cram\'er--Lundberg condition,
with $b=EX_1e^{\ga X_1}$, we have
\[
\frac{\tau(u)}{u}\to b^{-1} \qquad\mbox{in } P^{(u)}
\mbox{ probability }
\]
and
\[
\sup_{t\in[0,1]}\biggl\llvert \frac{X(t\tau(u))}{\tau(u)}-bt\biggr\rrvert
\to0 \qquad\mbox{in } P^{(u)}\mbox{ probability},
\]
indicating that ruin occurs due to the build up of small claims
which cause $X$ to behave as though it had positive drift; see Theorem
8.3.5 of \cite{EKM}.\hskip.2pt\footnote{The result cited in \cite{EKM} follows
from the work of Asmussen \cite{A1}, which is for the compound Poisson
model, but the result remains true for the general model.} By contrast
in the convolution equivalent case,
asymptotically,
ruin occurs in finite time (in distribution), and for ruin to occur,
the process must take a large jump from a neighborhood
of the origin to a neighborhood of $u$. This jump may result
in ruin, but if not, the resulting process $X-u$ subsequently behaves like
$X$ conditioned to hit $(0,\infty)$.
This representation of the limiting conditioned process leads to a
straightforward proof of (\ref{oSa}); see \cite{GM2}. However, the
description in the Cram\'er--Lundberg case is not sufficiently precise
to yield (\ref{oCr}).
What is needed is a more refined characterization of the process as
ruin approaches, specifically, a limiting description of the path from
the time of the last strict maximum before time $\tau(u)$ up until
time $\tau(u)$.

In the discrete time setting, such a result was proved by Asmussen
\cite{A1}.
Let $Z_k$ be i.i.d., nonlattice and set $S_n=Z_1+\cdots+Z_n$. Assume
the Cram\'er--Lundberg condition,
\[
\label{Crrw} Ee^{\ga Z_1}=1\quad\mbox{and}\quad EZ_1e^{\ga Z_1}<
\infty \qquad \mbox {for some } \ga>0.
\]
As above, let $\tau(u)$ be the first passage time of $S_n$ over level $u$
and $\gs(u)$ the time of the last strict ladder epoch prior to passage
[thus $\gs(0)=0$]. Set
\[
Z(u)=(Z_{\gs(u)+1},\dots, Z_{\tau(u)}).
\]
It follows from Section~8 of \cite{A1} that for $G$ bounded and continuous
%
\begin{eqnarray}
\label{As} && E^{(u)}G\bigl(Z(u), S_{\tau(u)}-u\bigr)
\nonumber\\[-8pt]\\[-8pt]\nonumber
&&\qquad \to\frac{\int_0^\infty e^{\ga y}E\{G(Z(0), S_{\tau(0)}-y); S_{\tau
(0)}>y, \tau(0)<\infty\} \,dy}{CE(S_{\tau(0)}e^{\ga S_{\tau(0)}};\tau(0)
<\infty)},
\end{eqnarray}
where $C=\lim_{u\to\infty}e^{\ga u}P(\tau(u)<\infty)$ and
$E^{(u)}$ denotes
expectation with respect to the conditional probability $P^{(u)}( \cdot )=P( \cdot \mid\tau(u)<\infty)$.
This result describes the limit of the conditioned process from the
time of the last strict ladder epoch prior to first passage over a high
level, up until the time of first passage. From it, the limiting
distribution of several quantities related to first passage, such as
those in~(\ref{oCr}), may be found in the random walk setting.

As it stands, the formulation in (\ref{As}) makes no sense for a
general L\'evy process. To apply even to the compound Poisson model,
the most popular risk model, some reformulation is needed. Furthermore,
to prove (\ref{As}), Asmussen derives a renewal equation by considering
the two cases $\tau(0)=\tau(u)$ and $\tau(0)<\tau(u)$. This is a standard
renewal theoretic device which has no hope of success in the general L\'
evy insurance risk model since typically $\tau(0)=0$. To circumvent
these problems, we apply arguments from fluctuation theory and
excursion theory. This allows us to describe, for any L\'evy process,
the final segment of the path from the time of the last maximum prior
to ruin, up until the time of ruin. This description is in terms of the
renewal measure $V$ of the ascending ladder height process and the
excursion measure of $X$ below its running supremum $\Xbar$. The key
observation that
ties together the two cases (\ref{Cr}) and (\ref{Sa}), and allows proof
of convergence as $u\to\infty$, is that in either case, $(\Vbar\circ
\ln) $ is regularly varying at infinity with index $-\ga$, where
$\Vbar
(u)=V(\infty)-V(u)$.
This allows us to derive not only new results in the Cram\'er--Lundberg
setting, but also to provide
a tool for transferring results from one setting to the other, and in
particular, to explain the striking similarity between results under
(\ref{Cr}) and (\ref{Sa}).

A very different description of the sample paths which lead to ruin
under (\ref{Cr}), can be found in Barczy and Bertoin \cite{BB}.
Building on results from Bertoin and Savov \cite{BS}, they
describe the sample paths in reverse time, from the time of ruin, in
terms of the associated exponentially tilted process conditioned to
stay positive and started with the limiting distribution of the
undershoot $u-X_{{\tau(u)}-}$. These two approaches are quite distinct and
the aim of \cite{BB} is somewhat different from here.
An interesting example related to the post ruin process is discussed in
\cite{BB}, but the paper is not specifically directed at insurance
risk. The limiting process here is described in forward time, and the
convergence is stronger than in \cite{BB}, in that it also applies to
certain discontinuous and unbounded functionals of the path.
Additionally, the results of \cite{BB} do not apply to the convolution
equivalent setting and so cannot explain the connection between results
such as (\ref{oCr}) and (\ref{oSa}).
The approach
in this paper may also prove useful in establishing
similar connections
for related processes.
For example, Mijatovic and Pistorius \cite{MP} recently showed that the
joint limit law of the undershoot and overshoot for the reflected
process under (\ref{Cr}) is the same as for the processes itself.
It now seems reasonable to conjecture that the analogous result holds
under (\ref{Sa}) and, furthermore, that this is a consequence of
a more general result related to
the sample path behavior of the process and the reflected processes under
(\ref{Cr}) and (\ref{Sa}) as first passage approaches.

Although not directly related to the current work, a description of the
sample paths leading to ruin has also been obtained for the heavy
tailed subexponential class of general L\'evy insurance risk processes.
This class was studied in the compound Poisson model by Asmussen and
Kl\"{u}ppelberg \cite{AK} and later for spectrally positive process by
Kl\"{u}ppelberg and Kyprianou \cite{KK}. Results for the general L\'evy
insurance risk process were obtained recently by Doney, Kl\"{u}ppelberg
and Maller \cite{DKM}.
The behavior of the paths is diametrically opposite to that in the
Cram\'er--Lundberg case, with ruin being a consequence of one extremely
large jump.

We conclude the \hyperref[sec1]{Introduction} with a brief outline of the paper.
Section~\ref{s2} contains the necessary fluctuation theory and
excursion theory to give a precise statement of the results. The main
results can then be found in Section~\ref{s4} together with an outline
of the general approach to their proof. Further results and proofs
related to Section~\ref{s2} are given in Section~\ref{s3} and the proof
of a preliminary result from Section~\ref{s4} is in Section~\ref{s10}.
Proof of the main convergence result under the Cram\'er--Lundberg
condition is given in Section~\ref{s5} and under the convolution
equivalent condition in Section~\ref{s6}. The special case where
$(0,\infty)$ is irregular is then briefly discussed in Section~\ref
{s7}. Specific calculations of limiting distributions as well as a
Gerber--Shiu EDPF are given in Section~\ref{s8}. Finally, the \hyperref[append]{Appendix}
contains a result in the case that $X$ is compound Poisson which, as is
often the case, needs to be treated separately.
Throughout $C, C_1, C_2,\dots$ will denote constants whose value is
unimportant and may change from one usage to the next.


\section{Fluctuation variables and excursion measure}\label{s2}

Let $L_t, {t\ge0}$, denote the local time at $0$ of the process
$X-\Xbar$, normalized by
%
\begin{equation}
\label{nLT} E\int_0^\infty
e^{-t}dL_t=1.
\end{equation}
Here, we are following Chaumont \cite{Ch} in our choice of
normalization. When $0$ is regular for $[0,\infty)$, $L$ is the unique
increasing, continuous, additive functional satisfying (\ref{nLT}) such
that the support of the measure $dL_t$ is the closure of the set $\{
t\dvtx \Xbar_t=X_t\}$ and $L_0=0$ a.s. If $0$ is irregular for $[0,\infty)$,
the set $\{s\dvtx  X_s>\overline{X}_{s-}\}$ of times of strict new maxima of
$X$ is discrete. Let $R_t=\llvert  \{s\in(0,t]\dvtx  X_s>\overline{X}_{s-}\}\rrvert  $ and
define the local time of $X-\Xbar$ at $0$ by
%
\begin{equation}
\label{LTi} L_t=\sum_{k=0}^{R_t}
\rme_k,
\end{equation}
where $\rme_k$, $k=0,1,\dots$ is an independent sequence of i.i.d.
exponentially distributed random variables with parameter
%
\begin{equation}
\label{nLT2} p=\frac{1}{1-E(e^{-\tau(0)};\tau(0)<\infty)}.
\end{equation}
Note that in this latter case, $dL_t$ has an atom of mass $\rme_0$ at
$t=0$ and thus the choice of $p$ ensures that (\ref{nLT}) holds.
Let $L^{-1}$ be the right continuous inverse of~$L$ and $H_s=\Xbar
_{L^{-1}_s}$. Then $(L^{-1}_s,H_s)_{s \geq0}$ is
the (weakly) ascending bivariate ladder process.

We will also need to consider the strictly ascending bivariate ladder
process, which requires a slightly different definition for $L$.
Specifically, when $0$ is regular for $(0,\infty)$, $L$ is the unique
increasing, continuous, additive functional as above.
When $0$ is irregular for $(0,\infty)$, $L$ is defined by (\ref{LTi}).
Thus, the only difference is for the compound Poisson process, where
the $L$ switches from being continuous to being given by (\ref{LTi}).
In this case, \textit{that is}, when $X$ is compound Poisson, the
normalization (\ref{nLT}) still holds, but now the support of the
measure $dL_t$ is the set of times of strict maxima of $X$, as opposed
to the closure of the set $\{t\dvtx \Xbar_t=X_t\}$.
$L^{-1}$ and $H$ are then defined as before in terms of $L$ and $\Xbar
$, and $(L^{-1}_s,H_s)_{s \geq0}$ is
the strictly ascending bivariate ladder process. See \cite{bert,doneystf} and particularly Chapter~6 of \cite{kypbook}.

In the following paragraph, $(L^{-1}_s,H_s)_{s \geq0}$ can be either
the weakly ascending or strictly ascending bivariate ladder
process.
When $X_t\to-\infty$ a.s., $L_\infty$ has an exponential distribution
with some parameter $q>0$, and the defective process $(L^{-1},H)$ may
be obtained from a nondefective process
$({\mathcal L}^{-1},{\mathcal H})$ by independent exponential killing at
rate $q > 0$.
We denote the bivariate L\'{e}vy measure of
$({\mathcal L}^{-1},{\mathcal H})$
by $\Pi_{{L}^{-1},{H}}(\cdot,\cdot)$.
The Laplace exponent $\kappa(a,b)$
of $(L^{-1},H)$,
defined by
\[
\label{kapdef} e^{-\kappa(a,b)} = E \bigl(e^{-a{L}^{-1}_1 -b{H}_1}; 1<L_\infty
\bigr) = e^{-q}E e^{-a{\mathcal L}^{-1}_1 -b{\mathcal H}_1}
\]
for values of $a,b\in\R$ for which the expectation is finite,
may be written
\[
\label{kapexp} \kappa(a,b) = q+\rmd_{L^{-1}}a+\rmd_Hb+\int
_{t\ge0}\int_{x\ge0} \bigl(1-e^{-at-bx}
\bigr) \Pi_{{L}^{-1}, {H}}(d t, d x),
\]
where $\rmd_{L^{-1}}\ge0$ and $\rmd_H\ge0$ are drift
constants. Observe that the normalization (\ref{nLT}) results in
$\kappa(1,0)=1$.
The bivariate renewal function of
$(L^{-1},H)$, given by
\begin{eqnarray*}
\label{Vkdef} V(t,x)&=& \int_0^\infty
e^{-qs}P\bigl({\mathcal L}_s^{-1}\le t,{\mathcal
H}_s\le x\bigr)\,d s
\\
&=& \int_0^\infty P\bigl({L}_s^{-1}
\le t,{H}_s\le x; s<L_\infty\bigr)\,d s,
\end{eqnarray*}
has Laplace transform
%
\begin{eqnarray}
\label{Vkap} \int_{t\ge0}\int_{x\ge0}e^{-at-bx}
V(d t,d x)&=&\int_{s= 0}^\infty e^{-qs}E
\bigl(e^{-a {\mathcal L}_s^{-1}-b{\mathcal H}_s}\bigr)\,d s
\nonumber\\[-8pt]\\[-8pt]\nonumber
&=& \frac
{1}{\kappa(a,b)},
\end{eqnarray}
provided $\gk(a,b)>0$. We will also frequently consider the renewal
function of $H$, defined on $\R$
by
\[
\label{VHdef} V(x)= \int_0^\infty
e^{-qs} P({\mathcal H}_s\le x)\,d s =\lim_{t\to\infty}V(t,x).
\]
Observe that $V(x)=0$ for $x<0$, while $V(0)>0$ iff ${\mathcal H}$ is
compound Poisson. Also
%
\begin{equation}
\label{Vinf} V(\infty):=\lim_{x\to\infty} V(x)=q^{-1}.
\end{equation}

From this point on, we will take $(L^{-1}, H)$ to be the strictly
ascending bivariate
ladder processes of $X$. Let $\whX_t=-X_t$, $t\ge0$ denote the dual
process, and $(\whL^{-1}, \whH)$ the weakly ascending bivariate
ladder processes of $\whX$. This is opposite to the usual convention,
and means some care needs to be taken when citing the literature in the
compound Poisson case.
This choice is made because it leads to more natural results and a
direct analogue of (\ref{As}) when $X$ is compound Poisson.
All quantities relating to $\whX$ will be denoted in the obvious way,
for example,
$\widehat{\tau}(0), \widehat{p}, \Pi_{{\whL}^{-1}, {\whH}}$,
$\widehat\kappa$
and $\whV$.
With these choices of bivariate ladder processes, together with the
normalization of the local times implying $\gk(1,0)=\whk(1,0)=1$,
the Wiener--Hopf factorization takes the form
%
\begin{equation}
\label{WH} {\kappa(a,-ib) \widehat\kappa(a,ib) = a-\Psi_X( b),\qquad
a\ge0, b\in\R.}
\end{equation}
If $\ga>0$ and $Ee^{\ga X_1}<\infty$, then by analytically extending
$\gk$, $\whk$ and $\Psi_X$, it follows from (\ref{WH}) that
\[
\label{WH1} {\kappa(a,-z) \widehat\kappa(a,z) = a-\Psi_X(-\rmi
z) \qquad\mbox{for } a\ge0, 0\le\Re z\le\ga.}
\]
If further $Ee^{\ga X_1}<1$, for example, when (\ref{Sa}) holds, then
$\Psi_X(-\rmi\ga)<0$ and since trivially $\whk(a,\ga)>0$, we have
%
\begin{equation}
\label{kpos} \gk(a,-\ga)>0 \qquad\mbox{for $a\ge0$.}
\end{equation}


Let $D$ be the Skorohod space of functions $w\dvtx [0,\infty)\to R$ which are
right continuous with left limits, equipped with the usual Skorohod
topology. The lifetime of a path $w\in D$ is defined to be $\zeta
(w)=\inf\{t\ge0\dvtx  w(s)=w(t)$ for all $s\ge t\}$, where we adopt the
standard convention that $\inf\varnothing=\infty$.
If $\gz(w)=\infty$ then $w(\gz)$ is taken to be some cemetery point.
Thus, for example, if $w(\gz)>y$ for some $y$ then necessarily $\gz
<\infty$. The jump in $w$ at time $t$ is given by $\Delta w_t=w(t)-w(t-)$.
We assume that $X$ is given as the coordinate process on $D$, and
the usual right continuous completion of the filtration generated by
the coordinate maps will be denoted $\{{\mathcal F}_t\}_{t\ge0}$.
$P_z$ is the probability measure induced on
${\mathcal F}=\bigvee_{t\ge0} {\mathcal F}_t$
by the L\'evy process starting at $z\in\R$, and
we usually write $P$ for $P_0$.\vadjust{\goodbreak}

Let\vspace*{1pt} $G=\{L^{-1}_{t-}\dvtx \Delta L^{-1}_{t}>0\}$ and $D=\{L^{-1}_{t}\dvtx \Delta
L^{-1}_{t}>0\}$ denote the set of left and right endpoints of excursion
intervals of $X-\Xbar$.
For $g\in G$, let $d\in D$ be the corresponding right endpoint of the
excursion interval ($d=\infty$ if the excursion has infinite
lifetime), and set
\[
\gep_g(t)=X_{(g+t)\wedge d}-\Xbar_g,\qquad t\ge0.
\]
Note, these are $X$-excursions in the terminology of Greenwood and
Pitman; see Remark 4.6 of \cite{GP}, as opposed to $X-\Xbar$ excursions.
Let
\[
\label{exc} {\mathcal E}=\bigl\{w\in D\dvtx  w(t)\le0 \mbox{ for all } 0\le t<
\zeta(w)\bigr\}
\]
and $\mathcal F^E$ the restriction of $\mathcal F$ to $\mathcal E$.
Then $\gep_g\in{\mathcal E}$ for each $g\in G$, and $\zeta(\gep_g)=d-g$.
The characteristic measure on $({\mathcal E}, {\mathcal F^E}) $ of the
$X$-excursions will be denoted $n$.

For fixed $u>0$, let
\[
G_{\tau(u)-}= \cases{ g, &\quad if $\tau(u)=d$ for some excursion interval
$(g,d)$
\cr
\tau(u), &\quad else.}
\]
If $X$ is compound Poisson, then $G_{\tau(u)-}$ is the first time of the
last maximum prior to $\tau(u)$. When $X$ is not compound Poisson,
$G_{\tau(u)-}$ is the left limit at $\tau(u)$ of $G_t=\sup\{s\le t\dvtx
X_s=\Xbar
_s\}$, explaining the reason behind this common notation.

Set
\[
Y_u(t) = X_{(G_{\tau(u)-}+t)\wedge\tau(u)} - \Xbar_{{\tau(u)-}}, \qquad t\ge0.
\]
Clearly, $\gz(Y_u)=\tau(u)-G_{\tau(u)-}$. If $\gz(Y_u)>0$ then
$G_{\tau(u)-}\in
G$, $\Xbar_{{\tau(u)-}}=\Xbar_{G_{\tau(u)-}}$ and $Y_u\in{\mathcal
E}$. If in
addition $\tau(u)<\infty$, equivalently $\gz(Y_u)<\infty$, then
$Y_u$ is
the excursion which leads to first passage over level $u$.
To cover the possibility that first passage does not occur at the end
of an excursion interval,
introduce
\[
\ceb={\mathcal E}\cup\{{\mathbf x}\dvtx  x\ge0\},
\]
where ${\mathbf x}\in D$ is the path which is identically $x$.
On the event $\gz(Y_u)=0$, that is $G_{\tau(u)-}=\tau(u)$, either
$X$ creeps
over $u$ in which case $Y_u={\mathbf0}$, or $X$ jumps over $u$ from its
current strict maximum in which case $Y_u={\mathbf x}$ where $x=\Delta
X_{\tau(u)}>0$ is the size of the jump at time $\tau(u)$. In all
cases, $Y_u\in
\ceb$.

Let\vspace*{1pt} ${\mathcal F^{\Ebar}}$ be the restriction of $\mathcal F$ to $\ceb
$. We extend $n$ trivially to a measure on ${\mathcal F^{\Ebar}}$ by
setting $n(\ceb\setminus{\mathcal E})=0$. Let $\tn$ denote the measure
on ${\mathcal F^{\Ebar}}$ obtained by pushing forward the measure $\Pi
_X^+$ with the mapping $x\to{\mathbf x}$, where $\Pi_X^+$ is the
restriction of $\Pi_X$ to $[0,\infty)$. Thus, $\tn({\mathcal E})=0$,
and for any Borel set $B\subset[0,\infty)$, $\tn(\{{\mathbf x}\dvtx x\in B\}
)=\Pi_X^+(B)$. Finally, let $\nbar=n+\rmd_{L^{-1}}\tn$.

For $u>0, s\ge0, y\ge0, \gep\in\ceb$
define
%
\begin{eqnarray}\label{Qu}
&& Q_u(ds, dy, d\gep)
\nonumber\\[-8pt]\\[-8pt]\nonumber
&&\qquad  = P\bigl(G_{\tau(u)-}\in ds,
u-\Xbar_{\tau(u)-}\in dy, Y_u\in d\gep, \tau(u)<\infty\bigr).
\end{eqnarray}
The starting point for our investigation is the following result, to be
proved in Section~\ref{s3}, which
provides a description of the sample paths from the (first) time of the
last maximum prior to\vadjust{\goodbreak} $\tau(u)$ until the time of first passage over
$u$. It
may be viewed as an extension of the quintuple law of Doney and
Kyprianou \cite{DK}; see the discussion following Proposition \ref{nmarg}.

\begin{teo}\label{ThmQ} For $u>0, s\ge0, y\ge0, \gep\in\ceb$,
%
\begin{eqnarray}
\label{teo1} Q_u(ds, dy, d\gep) &=& I(y\le u)V(ds, u-dy) \nbar
\bigl(d\gep,\gep(\gz )>y\bigr)
\nonumber\\[-8pt]\\[-8pt]\nonumber
&&{} +
\rmd_H \frac{\partial_-}{\partial_- u}V(ds,u)\delta_{0}(dy)\delta
_{{\mathbf0}}(d\gep),
\end{eqnarray}
where $\frac{\partial_-}{\partial_- u}$ denotes left derivative and
$\frac{\partial_-}{\partial_- u}V(ds,u)$ is the Lebesgue--Stieltges
measure associated with the function $\frac{\partial_-}{\partial_-
u}V(s,u)$ (which is increasing in s by (1.2) and~(3.5) of \cite{GM1}).
\end{teo}


\section{Statement of results and a unified approach}\label{s4}

In this section, we state the main results and outline a unified
approach to proving them under (\ref{Cr}) and (\ref{Sa}).
We assume from now on that $X_t\to-\infty$. We will be interested in a
marginalized version of (\ref{Qu}) conditional on $\tau(u)<\infty$. Thus,
for $u>0, y\ge0$ and $\gep\in\ceb$ define
\[
Q^{(u)}(dy, d\gep) = P^{(u)}(u-\Xbar_{\tau(u)-}\in dy,
Y_u\in d\gep),
\]
where recall $P^{(u)}( \cdot )=P( \cdot \mid\tau(u)<\infty)$.
Setting $\Vbar(u)=V(\infty)-V(u)$, and using the Pollacek--Khintchine formula,
%
\begin{equation}
\label{KP} P\bigl(\tau(u)<\infty\bigr)=q\Vbar(u),
\end{equation}
see Proposition 2.5 of \cite{kkm},
it follows from (\ref{teo1}) that
\[
 Q^{(u)}(dy, d\gep) = I(y\le u) \frac{V(u-dy)}{q\Vbar(u)} \nbar
\bigl(d\gep, \gep(\gz)>y\bigr) +\rmd_H \frac{V'(u)}{q\Vbar(u)}
\delta_{0}(dy)\delta _{{\mathbf
0}}(d\gep).
\]
Here, we have used the fact that $V$ is differentiable when $\rmd_H>0$,
see Theorem~VI.19 of \cite{bert}.
Now under either the Cram\'er--Lundberg condition (\ref{Cr}) or the
convolution equivalent condition (\ref{Sa}),
%
\begin{eqnarray}\label{vcon}
\frac{V(u-dy)}{q\Vbar(u)} &\stackrel{\mathrm{v}} {\longrightarrow }&
\frac{\ga}q e^{\ga y}\,dy\quad\mbox {and}\quad\rmd_H
\frac{ V'(u)}{q\Vbar(u)}\to\rmd_H\frac{\ga
}q
\nonumber\\[-12pt]\\[-8pt]
\eqntext{\mbox{as } u\to \infty;}
\end{eqnarray}
see Sections~\ref{s5} and \ref{s6} below. This suggests that under
suitable conditions on $F\dvtx [0,\infty)\times\ceb\to\R$,
%
\begin{equation}
\label{Fcon} \int_{[0,\infty)\times{\ceb} }F(y, \gep)Q^{(u)}(dy, d
\gep)\to \int_{[0,\infty)\times{\ceb} }F(y, \gep)Q^{(\infty)}(dy, d\gep),
\end{equation}
where
%
\begin{equation}
\label{Qinf} Q^{(\infty)}(dy, d\gep) = \frac{\ga}q
e^{\ga y}\,dy\, \nbar\bigl(d\gep, \gep(\gz)>y\bigr)+ \rmd_H
\frac{\ga}q \delta_{0}(dy)\delta_{{\mathbf0}}(d\gep),
\end{equation}
thus yielding a limiting description of the process as ruin approaches.
Observe that~(\ref{Fcon}) may be rewritten as
%
\begin{eqnarray}
\label{MR}
&& E^{(u)}F(u-\Xbar_{\tau(u)-},Y_u)
\nonumber\\[-8pt]\\[-8pt]\nonumber
&&\qquad \to\int
_{[0,\infty) }\frac{\ga}q e^{\ga y}\,dy \int
_{{\ceb}}F(y,\gep) \nbar\bigl(d\gep, \gep(\gz)>y\bigr)+
\rmd_H\frac{\ga}q F(0,{\mathbf0}),
\end{eqnarray}
indicating how the limiting behavior of many functionals of the process
related to ruin may be calculated.

To determine a broad class of functions $F$ for which (\ref{Fcon})
holds, first introduce
%
\begin{equation}
\label{h} h(y)=\int_{{\ceb}}F(y, \gep) \nbar\bigl(d\gep,
\gep(\gz)>y\bigr).
\end{equation}
We emphasize that throughout, $h$ will always depend on $F$, but we
will suppress this dependence to ease notation.
Since, by (\ref{vcon}), (\ref{Fcon}) is equivalent to
%
\begin{equation}
\label{conh} \int_0^u h(y)
\frac{V(u-dy)}{q\Vbar(u)}
\to\int_0^\infty h(y)\frac{\ga e^{\ga y} }{q}\,dy,
\end{equation}
it will be of interest to know when $h$ is continuous a.e.$ $with
respect to Lebesgue measure $m$.
The most obvious setting in which
the condition on $B_y$ below holds, is when $F$ is continuous
in $y$ for each $\gep$. In particular, it holds when $F$ is jointly
continuous. The boundedness condition holds when $F$ is bounded, but
applies to certain unbounded functions.

\begin{prop}\label{hcont} Assume $Ee^{\ga X_1}<\infty$, and
$F\dvtx [0,\infty
)\times{\ceb}\to\R$ is product measurable, and $F(y, \gep)e^{-\ga
(\gep
(\gz)-y)}I(\gep(\gz)>y)$ is bounded in $(y,\gep)$. Further assume
$\nbar
(B_y^c)=0$ for a.e. y with respect to Lebesgue measure $m$, where
\[
B_y=\bigl\{\gep\dvtx  F(\cdot ,\gep) \mbox{ is continuous at } y\bigr
\}.
\]
Then $h$ is continuous a.e. w.r.t. $m$.
\end{prop}

We can now state the main results.

\begin{teo}\label{CrThm}
If (\ref{Cr}) holds, and $F\ge0$ satisfies the hypotheses for
Proposition \ref{hcont}, then (\ref{Fcon}), equivalently (\ref{MR}), holds.
In particular,
\[
\label{Qconw} Q^{(u)}(dy, d\gep)\stackrel{\mathrm{w}} {
\longrightarrow}Q^{(\infty
)}(dy, d\gep).
\]
\end{teo}

\begin{teo}\label{SaThm}
If (\ref{Sa}) holds, $F\ge0$ satisfies the hypotheses of Proposition~\ref{hcont}, and
%
\begin{equation}
\label{FSa} F(y,\gep)I\bigl(\gep(\gz)>y\bigr) \to0 \qquad\mbox{uniformly in }
\gep\in {\ceb }\mbox{ as } y\to\infty,
\end{equation}
then (\ref{Fcon}), equivalently (\ref{MR}), holds. Condition (\ref{FSa}) holds if, for example,
$F$ has compact support. In particular,
%
\begin{equation}
\label{Qconv} Q^{(u)}(dy, d\gep)\stackrel{\mathrm{v}} {
\longrightarrow}Q^{(\infty
)}(dy, d\gep).
\end{equation}
\end{teo}

While the outline of the proofs of these two results is the same, the
details need to be handled differently. In the Cram\'er--Lundberg case,
the key renewal theorem is used, whereas in the convolution equivalent
case, special properties of convolution equivalent distributions are
used. This results in different classes of functions for which (\ref
{Fcon}) holds.
The extra condition (\ref{FSa}) in Theorem \ref{SaThm} cannot be
dispensed with, else the convergence in (\ref{Qconv}) could be
improved to weak convergence. However, as we show later in (\ref
{Qmass}), the total mass of $Q^{(\infty)}$ under (\ref{Sa}) is less
than one.
The convergence in (\ref{MR}) may be expressed alternatively in terms
of the overshoot $X_{\tau(u)}-u$ rather than the undershoot $u-\Xbar
_{\tau(u)
-}$, making it more analogous to~(\ref{As}); see Theorems \ref{GCr} and
\ref{GSa}.

Evaluation, or even simplification, of the limit in (\ref{MR}) for a
specific functional of the path is, in general, difficult to achieve.
Here, we give an example where it is possible and which arises
naturally in risk theory. Note that the function $f$ below can grow
exponentially in the overshoot variable. This allows for the
calculation of certain unbounded Gerber--Shiu expected discounted
penalty functions;
see Section~\ref{s8}.

\begin{teo}\label{thmf} Assume that (\ref{Cr}) holds, and
$f\dvtx [0,\infty
)^4\to[0,\infty)$ is a Borel function which is jointly continuous in
the first three variables and $e^{-\ga x}f(y,x,v,t)$ is bounded. Then
%
\begin{eqnarray}
\label{thmf1} \quad&& E^{(u)}f\bigl( u-\Xbar_{\tau(u)-},X_{\tau(u)}-u,
u-X_{\tau(u)-}, \tau (u)- G_{\tau(u)-}\bigr)\nonumber
\\
&&\qquad \to
\int_{x\ge0}\int_{y\ge0}\int
_{v\ge0}\int_{t\ge0} f(y,x, v,t)
\nonumber\\[-8pt]\\[-8pt]\nonumber
&&\hspace*{122pt}{}\times
\frac{\ga}q e^{\ga y} \,dy\, I(v\ge y) \whV(dt,dv-y)
\Pi_X(v+dx)
\\
&&\quad\qquad{} +
\rmd
_H\frac{\ga}q f(0,0,0,0).\nonumber
\end{eqnarray}
In particular, we have joint convergence; for $y\ge0, x\ge0, v\ge0,
t\ge0$
%
\begin{eqnarray}
\label{jtcon} &&P^{(u)}\bigl(u-\Xbar_{\tau(u)-}\in dy,
X_{\tau(u)}-u\in dx,u-X_{\tau
(u)-}\in dv,\nonumber
\\
&&\hspace*{166pt}  \tau(u)- G_{\tau(u)-}
\in dt\bigr)
\nonumber\\[-8pt]\\[-8pt]\nonumber
&&\qquad \stackrel{\mathrm{w}} {\longrightarrow}\frac{\ga}q e^{\ga y} \,dy\,
I(v\ge y) \whV(dt,dv-y)\Pi_X(v+dx)
\\
&&\hspace*{9pt}\quad\qquad{}  + \rmd _H
\frac{\ga}q \delta_{(0,0,0,0)}(dx,dy,dv,dt).\nonumber
\end{eqnarray}
\end{teo}

In the convolution equivalent setting, the following result
extends Theorem 10 of \cite{DK}.

\begin{teo}\label{thmfSa} Assume (\ref{Sa}) holds and that
$f\dvtx [0,\infty
)^4\to[0,\infty)$ satisfies~(\ref{fun0}), is jointly continuous in the
first three variables, $e^{-\ga x}f(y,x,v,t)$ is bounded, and
\[
\sup_{x>0,t\ge0, v\ge y}f(y,x,v,t)\to0 \qquad\mbox{as } y\to
\infty.
\]
Then
%
\begin{eqnarray}
\label{Safcon} \quad&& E^{(u)}f\bigl(u-\Xbar_{\tau(u)-},X_{\tau(u)}-u,
u-X_{\tau(u)-}, \tau (u)- G_{\tau(u)-}\bigr)\nonumber
\\
&&\qquad \to
\int_{x\ge0}\int_{y\ge0}\int _{v\ge0}\int_{t\ge0} f(x,y,v,t)
\nonumber\\[-8pt]\\[-8pt]\nonumber
&&\hspace*{121pt}{}\times
\frac{\ga}q e^{\ga y} \,dy\, I(v\ge y) \whV(dt,dv-y)
\Pi_X(v+dx)
\\
&&\quad\qquad{}+\rmd
_H\frac{\ga}q f(0,0,0,0).\nonumber
\end{eqnarray}
In particular, we have joint convergence; for $y\ge0, x\ge0, v\ge0,
t\ge0$
%
\begin{eqnarray}
\label{jtconv}
&&P^{(u)}\bigl(u-\Xbar_{\tau(u)-}\in dy, X_{\tau(u)}-u\in dx, u-X_{\tau
(u)-}\in dv,\nonumber
\\
&&\hspace*{166pt} \tau(u)- G_{\tau(u)-} \in dt\bigr)
\nonumber\\[-8pt]\\[-8pt]\nonumber
&&\qquad \stackrel{\mathrm{v}} {\longrightarrow}\frac{\ga}q e^{\ga y} \,dy\,
I(v\ge y) \whV(dt,dv-y)\Pi_X(v+dx)
\\
&&\hspace*{9pt}\quad\qquad{}  + \rmd _H
\frac{\ga}q \delta_{(0,0,0,0)}(dx,dy,dv,dt).\nonumber
\end{eqnarray}
\end{teo}

Theorems \ref{CrThm} and \ref{SaThm} describe, in a very general sense,
how to transfer results from the Cram\'er--Lundberg setting to the
convolution equivalent setting and vice versa. Theorems \ref{thmf} and
\ref{thmfSa} provide a specific example of this. However, since the
mode of convergence is $\stackrel{\mathrm{w}}{\longrightarrow}$
under (\ref{Cr}) and $\stackrel{\mathrm{v}}{\longrightarrow}$ under
(\ref{Sa}), some subtleties may arise. For example,
the marginal distributions of the limit in (\ref{jtcon}) can be readily
calculated using (\ref{DKV}) below, and
consequently under (\ref{Cr}) we obtain, in addition to~(\ref
{oCr}),\footnote{Strictly speaking, the proof of (\ref{oCr}) in \cite
{GMS} assumes that $(L^{-1}, H)$ is the weakly ascending ladder
process, whereas the marginals of (\ref{jtcon}) yield the same formulae
as (\ref{oCr}) but with $(L^{-1}, H)$ the strictly ascending ladder
process. Thus, as can be easily checked directly, the limiting
expressions must agree irrespective of the choice of ascending ladder process.
This remark applies to (\ref{oSa}) and several other limiting
distributions discussed here.}
\begin{eqnarray*}
\label{uCr} %
P^{(u)}(u-X_ {{\tau(u)}-}
\in dx)&\stackrel{\mathrm {w}} {\longrightarrow} &q^{-1} \alpha
\rmd_H\delta_{0}(dx)
\\
&&{} +q^{-1} \alpha
e^{\alpha x} \pibar_X(x)\,dx\int_{0\le v\le x}e^{-\ga v}
\whV(dv),
\\
 P^{(u)}\bigl(\tau(u)-G_ {{\tau(u)}-}\in dt\bigr)
&\stackrel {\mathrm{w}} {\longrightarrow}& q^{-1} \bigl(\alpha\rmd
_H\delta_0(d t)+K(dt) \bigr),
\end{eqnarray*}
where
%
\begin{equation}
\label{K} K(dt)=\int_{z\ge0}\bigl(e^{\alpha z} -1
\bigr)\Pi_{L^{-1},H}(dt, dz).
\end{equation}
Under (\ref{Sa}), some care is needed. The marginals of the limit in
(\ref{jtconv}) are the same as in (\ref{jtcon}), but they all have mass
less than one. This does not mean that we can simply replace weak
convergence of the marginals under (\ref{Cr}) with vague convergence
under (\ref{Sa}).
For the undershoots of $X$ and $\Xbar$, this is correct, but the
overshoot and $\tau(u)-G_ {{\tau(u)}-}$ both converge weakly under
(\ref{Sa});
indeed they converge jointly, as will be shown in Proposition
\ref{p6}.
Consequently, an extra term appears in the limit of the overshoot in
(\ref{oSa}) to account for the missing mass. Similarly, for $\tau(u)-G_
{{\tau(u)}-}$, see
(\ref{uSa}).

Based on the outline of the proofs of Theorems \ref{CrThm} and \ref
{SaThm} given at the beginning of this section, it is natural to ask if
any other limits are possible in (\ref{vcon}), thus leading to
different forms of the limit in (\ref{MR}). However, this is not the case.
More precisely, if ${V(u-dy)}/{\Vbar(u)}$ converges vaguely to a
nonzero locally finite Borel measure, then ${\Vbar(u-y)}/{\Vbar(u)}$
converges as $u\to\infty$ on a dense set of $y$. Hence, by Theorem
1.4.3 of \cite{BGT}, $\Vbar(\ln u)$ is regularly varying at infinity
with some index $-\ga$. Consequently, the limit in (\ref{vcon}) must be
of the form given. The only general classes of processes that the
author is aware of which satisfy (\ref{vcon}) are those studied in this
paper, namely the Cram\'er--Lundberg and convolution equivalent cases.


\section{Proof of Theorem \texorpdfstring{\protect\ref{ThmQ}}{2.1} and related results}\label{s3}

The following result will be needed in the proof of Theorem \ref{ThmQ}.

\begin{prop}\label{Ch1} If $X$ is not compound Poisson, then for $s\ge
0, x\ge0$,
\[
\label{ch1} \rd_{L^{-1}}V(ds, dx)=P(\Xbar_s=X_s
\in dx)\,ds.
\]
\end{prop}

\begin{pf}   For any $s\ge0, x\ge0$
\begin{eqnarray*}
\label{ch1a} \rd_{L^{-1}} \int_0^{L_\infty} I
\bigl(L^{-1}_t\le s, H_t\le x\bigr) \,dt &=&
\rd_{L^{-1}} \int_0^{L_\infty} I
\bigl(L^{-1}_t\le s, \Xbar_{L^{-1}_t}\le x\bigr) \,dt
\\
&=& \rd_{L^{-1}} \int_0^s I(
\Xbar_r\le x) \,dL_{r}
\\
&=& \int_0^s I(\Xbar_r\le x)I(
\Xbar_r=X_r) \,dr,
\end{eqnarray*}
by Theorem 6.8 and Corollary 6.11 of\vspace*{1pt} Kyprianou \cite{kypbook}, which
apply since $X$ is not compound Poisson, [Kyprianou's $(L^{-1},H)$ is
the weakly ascending ladder process in which case the result holds in
the compound Poisson case also]. Taking expectations completes the proof.
\end{pf}

\begin{pf*}{Proof of Theorem \ref{ThmQ}}
There are three possible
ways in which $X$ can first cross level $u$; by a jump at the end of an
excursion interval, by a jump from a current strict maximum or by
creeping. We consider each in turn.

Let $f, h$ and $j$ be nonnegative bounded continuous functions. Since
$\Xbar_{L^{-1}_{t-}}$ is left continuous, we may apply the
master formula of excursion theory, Corollary IV.11 of \cite{bert}, to obtain
%
\begin{eqnarray}
\label{t1} && E\bigl\{f(G_{\tau(u)-})h(u-\Xbar_{\tau(u)-})j(Y_u);X_{\tau(u)}>u,
G_{\tau(u)-}<\tau(u)<\infty\bigr\}\nonumber
\\
&&\qquad = E\sum_{g\in G}f(g)h(u-\Xbar_{g})j(
\gep_g)I\bigl(\Xbar_{g}\le u, \gep _g(
\zeta )>u-\Xbar_{g}\bigr)\nonumber
\\
&&\qquad =E\int_0^\infty dL_t\int
_{\mathcal E} f(t)h(u-\Xbar_{t})j(\gep )I\bigl(\Xbar
_{t}\le u, \gep(\zeta)>u-\Xbar_{t}\bigr)n(d\gep)\nonumber
\\
&&\qquad =\int_{\mathcal E} j(\gep)n(d\gep) E\int_0^{L_\infty
}f
\bigl(L^{-1}_r\bigr)h(u-H_r)I
\bigl(H_r\le u, \gep(\zeta)>u-H_r\bigr) \,dr\hspace*{-15pt}
\\
&&\qquad =\int_{\mathcal E} j(\gep)n(d\gep) \int_{s\ge0}
\int_{0\le y\le u} f(s)h(u-y)I\bigl(\gep(\zeta)>u-y\bigr) V(ds, dy)\nonumber
\\
&&\qquad =\int_{\mathcal E} j(\gep)n(d\gep) \int_{s\ge0}
\int_{0\le y\le u} f(s)h(y)I\bigl(\gep(\zeta)>y\bigr) V(ds, u-dy)\nonumber
\\
&&\qquad =\int_{s\ge0}\int_{0\le y\le u}\int
_{\mathcal E} f(s)h(y)j(\gep )V(ds, u-dy)n\bigl(d\gep, \gep(\zeta)>y
\bigr).\nonumber
\end{eqnarray}

Next, define $\tj\dvtx [0,\infty)\to\R$ by $\tj(x)=j({\mathbf x})$. Then, since
$Y_u(t)=\Delta X_{\tau(u)}$ for all $t\ge0$ on
$\{X_{\tau(u)}>u,G_{\tau(u)-}=\tau(u)<\infty\}$, we have by the
compensation formula,
%
\begin{eqnarray}
\label{t2} &&E\bigl\{f(G_{\tau(u)-})h(u-\Xbar_{\tau(u)-})j(Y_u);X_{\tau
(u)}>u,G_{\tau(u)-}=
\tau(u)<\infty\bigr\}\nonumber
\\
&&\qquad= E\bigl\{f(G_{\tau(u)-})h(u-\Xbar_{\tau(u)-})\tj(\Delta
X_{\tau
(u)});\nonumber
\\
&&\hspace*{70pt} X_{\tau(u) }>u,G_{\tau(u)-}=\tau(u)<\infty\bigr\}\nonumber
\\
&&\qquad=E\sum_{s}f(s)h(u-\Xbar_{s-})
\tj(\Delta X_s)I(X_{s-}=\Xbar _{s-}\le u, \Delta
X_s>u-\Xbar_{s-})
\nonumber\\[-8pt]\\[-8pt]\nonumber
&&\qquad=E\int_0^\infty f(s)h(u-
\Xbar_s)I(X_s=\Xbar_s\le u)\,ds\int
_{\xi
}\tj(\xi)I(\xi>u-\Xbar_s)
\Pi_X(d\xi)
\\
&&\qquad=\int_0^\infty f(s)\int
_{0\le y\le u} h(y)\int_{\xi}\tj (\xi )I(\xi>y)
\Pi_X(d\xi)P(\Xbar_s=X_s\in u- dy)\,ds\nonumber
\\
&&\qquad=\int_{s\ge0}f(s)\int_{0\le y\le u}h(y)\int
_{\ceb} j(\gep )\tn \bigl(d\gep, \gep(\zeta)>y\bigr)P(
\Xbar_s=X_s\in u- dy)\,ds\nonumber
\\
&&\qquad=\int_{s\ge0}\int_{0\le y\le u}\int
_{\ceb} f(s)h(y)j(\gep )\tn \bigl(d\gep, \gep(\zeta)>y\bigr)
\rmd_{L^{-1}}V(ds, u-dy),\nonumber
\end{eqnarray}
where the final equality follows from Proposition \ref{Ch1} if $X$ is
not compound Poisson. If $X$ is compound Poisson the first and last
formulas of (\ref{t2}) are equal because $P(G_{\tau(u)-}=\tau(u))=0$
and $\rmd _{L^{-1}}=0$ [recall that $(L^{-1},H) $ is the strictly ascending
ladder process].

Finally,
%
\begin{eqnarray}
\label{t3} && E\bigl\{f(G_{\tau(u)-}) h(u-\Xbar_{\tau(u)-})j(Y_u);X_{\tau(u)}=u,
\tau (u)<\infty\bigr\}\nonumber
\\
&&\qquad = h(0)j({\mathbf0})E\bigl\{f\bigl({\tau(u)}\bigr);X_{\tau(u)}=u, \tau(u)<
\infty\bigr\}
\\
&&\qquad =\rmd_H h(0)j({\mathbf0})\int_s f(s)
\frac{\partial_-}{\partial_- u}V(ds,u)\nonumber
\end{eqnarray}
if $\rmd_H>0$, by (3.5) of \cite{GM1}. If $\rmd_H=0$, then $X$ does not
creep, and so $P(X_{\tau(u)}=u)=0$. Thus, (\ref{t3}) holds in this
case also.
Combining the three terms (\ref{t1}), (\ref{t2}) and (\ref{t3}) gives
the result.
\end{pf*}

The next two results will be used to calculate limits such as those of
the form (\ref{jtcon}) and (\ref{jtconv}).

\begin{prop}\label{2222} For $t\ge0, z\ge0$,
%
\begin{equation}
\label{ch2} \whV(dt, dz)=n\bigl(\gep(t)\in-dz, \zeta>t\bigr)\,dt+\rd_{L^{-1}}
\gd_{(0,0)}(dt,dz).
\end{equation}
\end{prop}

\begin{pf}   If $X$ is not compound Poisson nor $\llvert  X\rrvert  $ a
subordinator, (\ref{ch2}) follows from (5.9) of \cite{Ch} applied to
the dual process $\whX$.

If $X$ is a subordinator, but not compound Poisson, then $n$ is the
zero measure and $\rd_{L^{-1}}=1$ by (\ref{nLT}). On the other hand,
$({\whL}^{-1},{\whH})$ remains at $(0,0)$ for an exponential amount of
time with parameter $\widehat p=1$, by (\ref{nLT2}), and is then killed.
Hence, (\ref{ch2}) holds.

If $-X$ is a subordinator, then
$({\whL}_t^{-1},{\whH}_t)=(t,\whX_t)$ and so $\whV(dt, dz)=P(X_t\in
-dz)\,dt$. On the other hand, $n$ is proportional to the first, and only,
excursion, so $n(\gep(t)\in-dz, \zeta>t)=cP(X_t\in-dz)$ for some
$c>0$. Since $\rd_{L^{-1}}=0$, we thus only need check that $\llvert  n\rrvert  =1$. But
$G=\{0\}$, and so by the master formula
\[
1=E\sum_{g\in G} e^{-g}= E\int
_0^\infty e^{-t} \,dL_t\int
_{\mathcal
E}n(d\gep)=\llvert n\rrvert .
\]

To complete the proof, it thus remains to prove (\ref{ch2}) when $X$ is
compound Poisson. We defer this case to the \hyperref[append]{Appendix}.
\end{pf}

For notational convenience, we define $\gep(0-)=0$ for $\gep\in\ceb
$.
Thus, in particular, ${\mathbf x}(\gz-)=0$ since $\gz({\mathbf x})=0$. Note
also that ${\mathbf x}(\gz)=x$.

\begin{prop}\label{nmarg} For $t\ge0$, $z\ge0$ and $x>0$,
%
\begin{equation}
\label{nmarg1} \nbar\bigl(\gz\in dt, \gep(\gz-)\in-dz, \gep(\gz)\in dx\bigr)=
\whV (dt,dz)\Pi_X(z+dx).
\end{equation}
\end{prop}

\begin{pf}   First consider the case $t>0$, $z\ge0$ and
$x>0$. For any $0<s<t$, using the Markov property of the excursion
measure $n$, we have
\begin{eqnarray*}
&& \nbar\bigl(\gz\in dt, \gep(\gz-)\in-dz, \gep(\gz)\in dx\bigr)\nonumber
\\
&&\qquad = \int_{y\ge0}n\bigl(\gep(s)\in-dy, \zeta>s
\bigr)\nonumber
\\
&&\hspace*{52pt}{}\times P_{-y}\bigl({\tau(0)}\in dt-s, X_{\tau(0)-}\in-dz,
X_{\tau(0)}\in dx\bigr)
\\
&&\qquad = \int_{y\ge0} n\bigl(\gep(s)\in-dy, \zeta>s\bigr)\nonumber
\\
&&\hspace*{20pt}\quad\qquad{}\times P\bigl({\tau(y)}\in dt-s, X_{\tau(y)-}\in y-dz,
X_{\tau(y)}\in y+dx\bigr).\nonumber
\end{eqnarray*}
By the compensation formula, for any positive bounded Borel function $f$,
\begin{eqnarray*}
&& E\bigl\{f\bigl(\tau(y)\bigr);  X_{\tau(y)-}\in y-dz, X_{\tau(y)}\in
y+dx\bigr\}
\\
&&\qquad = E \biggl\{\sum_r f(r); \Xbar_{r-}
\le y, X_{r-}\in y-dz, X_r\in y+dx\biggr\}
\\
&&\qquad = \int_{r=0}^\infty f(r)P( \Xbar_{r-}
\le y, X_{r-}\in y-dz) \,dr\,\Pi _X(z+dx).
\end{eqnarray*}
Thus,
\begin{eqnarray*}
&& P\bigl({\tau(y)}\in dt-s, X_{\tau(y)-}\in y-dz, X_{\tau(y)}\in y+dx
\bigr)
\\
&&\qquad =P( \Xbar_{t-s}\le y, X_{t-s}\in y-dz) \,dt\, \Pi_X(z+dx).
\end{eqnarray*}
Hence,
\begin{eqnarray*}
&& \nbar \bigl(\gz\in dt, \gep(\gz-)\in-dz, \gep(\gz)\in dx\bigr)
\\
&&\qquad = \int_{y\ge0}n\bigl(\gep(s)\in-dy, \zeta>s\bigr)P(
\Xbar_{t-s}\le y, X_{t-s}\in y-dz) \,dt\, \Pi_X(z+dx)
\\
&&\qquad = n\bigl(\gep(t)\in-dz, \zeta>t\bigr)\,dt\, \Pi_X(z+dx)
\\
&&\qquad = \whV(dt,dz) \Pi_X(z+dx)
\end{eqnarray*}
by (\ref{ch2}).

Finally,
if $t=0$, then for any positive bounded Borel function,
\begin{eqnarray*}
&& \int_{\{(t,z,x)\dvtx t=0, z\ge0,x>0\}} f(t,z , x)\nbar\bigl(\gz\in dt, \gep (\gz -)
\in-dz, \gep(\gz)\in dx\bigr)
\\
&&\qquad = \rmd_{L^{-1}}\int_{x>0}f(0,0,x )\tn\bigl(\gep(\gz)
\in dx\bigr)
\\
&&\qquad = \rmd_{L^{-1}}\int_{x>0}f(0,0, x)
\Pi_X^+(dx)
\\
&&\qquad = \int_{\{(t,z,x)\dvtx t=0, z\ge0,x>0\}} f(t,z , x)\whV(dt,dz)\Pi_X(z+dx)
\end{eqnarray*}
by Proposition \ref{2222}.
\end{pf}

As mentioned earlier, Theorem \ref{ThmQ} may be viewed as an extension
of the quintuple law of Doney and Kyprianou \cite{DK}. To see this,
observe that from Theorem~\ref{ThmQ} and Proposition \ref{nmarg}, for
$u>0, s\ge0,t\ge0, 0\le y\le u\wedge z, x\ge0$,
%
\begin{eqnarray}
\label{QL}
&& P\bigl(G_{\tau(u)-} \in ds, u-\Xbar_{\tau(u)-}\in dy,
\tau(u)-G_{\tau
(u)-}\in dt,\nonumber
\\
&&\hspace*{95pt} u-X_{\tau(u)-}\in dz, X_{\tau
(u)}-u
\in dx\bigr)\nonumber
\\
&&\qquad = P\bigl(G_{\tau(u)-}\in ds, u-\Xbar_{\tau(u)-}\in dy,
\zeta(Y_u)\in dt,\nonumber
\\
&&\hspace*{84pt} Y_u(\zeta-)\in y-dz,
Y_u(\zeta )\in y+dx\bigr)
\nonumber\\[-8pt]\\[-8pt]\nonumber
&&\qquad = I(x>0)V(ds, u-dy) \nbar\bigl(\gz\in dt, \gep(\gz-)\in y-dz, \gep(\gz )\in
y+dx\bigr)
\\
&&\quad\qquad{} + \rmd_H \frac{\partial_-}{\partial_- u}V(ds,u)
\delta_{(0,0,0,0)}(dt,dx,dz,dy)\nonumber
\\
&&\qquad = I(x>0)V(ds, u-dy)\whV(dt,dz-y)\Pi_X(z+dx)\nonumber
\\
&&\quad\qquad{} +\rmd_H \frac{\partial_-}{\partial_- u}V(ds,u)
\delta_{(0,0,0,0)}(dt,dx,dz,dy).\nonumber
\end{eqnarray}
When $X$ is not compound Poisson, this is the statement of Theorem 3 of
\cite{DK} with the addition of the term due to creeping; see Theorem
3.2 of \cite{GM1}.
When $X$ is compound Poisson the quintuple law, though not explicitly
stated in \cite{DK}, remains true and can be found in \cite{EK}. In
that case, the result is slightly different from (\ref{QL}) since the
definitions of $G_{\tau(u)-}, V$ and $\whV$ then differ due to the choice
of $(L^{-1},H)$ as the weakly ascending ladder process in \cite{DK} and
\cite{EK}.
Thus, we point out that Vigon's \'equation amicale invers\'ee, \cite{Vig},
%
\begin{equation}
\label{EA} \Pi_H(dx)= \int_{z\ge0} \whV(dz)
\Pi_X(z+dx),\qquad x>0,
\end{equation}
and Doney and Kyprianou's extension,
%
\begin{equation}
\label{DKV} \Pi_{{L}^{-1},{H}}(dt,dx)= \int_{v\ge0}
\whV(dt, dv)\Pi_X(v+dx),\qquad x>0,t\ge0,
\end{equation}
continue to hold with our choice of $(L^{-1},H)$ as the strongly
ascending ladder process. The proof of (\ref{DKV}) is analogous to the
argument in Corollary 6 of \cite{DK}, using~(\ref{QL}) instead of Doney
and Kyprianou's quintuple law, and (\ref{EA}) follows immediately from
(\ref{DKV}).

\begin{cor}\label{nV}
For $x>0$,
%
\begin{equation}
\label{nV1} \nbar\bigl(\gep(\gz)\in dx\bigr)=\Pi_H(dx).
\end{equation}
\end{cor}

\begin{pf}   Integrating out in (\ref{nmarg1}),
\[
\label{Vig} \nbar\bigl(\gep(\gz)\in dx\bigr)=\int_{z\ge0}
\whV(dz)\Pi_X(z+dx)=\Pi_H(dx),
\]
by (\ref{EA}).
\end{pf}


\section{Proof of Proposition \texorpdfstring{\protect\ref{hcont}}{3.1}}\label{s10}

From (\ref{h}),
\begin{eqnarray*}
h(z)&=& \int_{{\ceb}}f_z(\gep)\nbar(d\gep),
\end{eqnarray*}
where
\[
f_z(\gep)=F(z, \gep)I\bigl(\gep(\gz)>z\bigr).
\]
Fix $y>0$ and assume $\llvert  z-y\rrvert  <y/2$. Then for some constant $C$,
independent of $z$ and $\gep$,
%
\begin{equation}
\label{h1} \bigl\llvert f_z(\gep)\bigr\rrvert \le C
e^{\ga(\gep(\gz)-z)}I\bigl(\gep(\gz)>z\bigr)\le C e^{\ga
(\gep(\gz
)-y/2)}I\bigl(\gep(
\gz)>y/2\bigr).
\end{equation}
By (\ref{nV1}),
%
\begin{equation}
\label{h2} \int_{{\ceb}}e^{\ga(\gep(\gz)-y/2)}I\bigl(\gep(\gz)>y/2
\bigr)\nbar(d\gep )=\int_{x>y/2}e^{\ga(x-y/2)}
\Pi_H(dx),
\end{equation}
and since $Ee^{\ga X_1}<\infty$, this last integral is finite by
Proposition 7.1 of \cite{G}.

Now let $A=\{y\dvtx \nbar(B_y^c)=0\}$ and $C_H=\{y\dvtx \Pi_H(\{y\})=0\}$. Then
$C_H^c$ is countable and
\[
\nbar\bigl(\gep(\gz)=y\bigr)=\Pi_H\bigl(\{y\}\bigr)=0\qquad
\mbox{if } y\in C_H.
\]
Thus, if $y>0$, $y\in A\cap C_H$ and $z\to y$, then $f_z(\gep)\to
f_y(\gep)$ a.e. $ \nbar$. Hence, by (\ref{h1}) and (\ref{h2}), we can
apply dominated convergence to obtain continuity of $h$ at such~$y$.
Since $m(A^c)=m(C_H^c)=0$, this completes the proof.


\section{Proofs under the Cram\'er--Lundberg condition}\label{s5}

In studying the process $X$ under the Cram\'er--Lundberg condition
(\ref{Cr}), it is useful to introduce the Esscher transform.
Thus, let $P^*$ be the measure on $\mathcal F$ defined by
\[
\label{P*tau} dP^*=e^{\alpha X_t} \,dP \qquad\mbox{on } {\mathcal
F}_{t},
\]
for all $t\ge0$. Then $X$ under $P^*$ is the Esscher transform of $X$.
It is itself a L\'evy process with $E^*X>0$; see Section~3.3 of
Kyprianou \cite{kypbook}.

When (\ref{Cr}) holds, Bertoin and Doney \cite{BD94} extended the
classical Cram\'er--Lundberg estimate for ruin to a general L\'evy
process; assume $X$ is nonlattice in the case that $X$ is compound
Poisson, then
%
\begin{equation}
\label{BD} \lim_{u\to\infty} e^{\ga u}P\bigl(\tau(u)<
\infty\bigr)=\frac{q}{\ga m^* },
\end{equation}
where $m^*=E^*H_1$.
Under $P^*$, $H$ is a nondefective subordinator with drift $\mathrm{d}^*_H$ and
L\'evy measure $\Pi^*_H$ given by
\[
\mathrm{d}^*_H=\mathrm{d}_H \quad\mbox{and}\quad\Pi^*_H(dx)=e^{\ga x}
\Pi_H(dx),
\]
and so
%
\begin{equation}
\label{m*} E^*H_1= \mathrm{d}^*_H + \int
_0^\infty x \Pi^*_H(dx)=
\mathrm{d}_H+\int_0^\infty x
e^{\ga
x} \Pi_H(dx).
\end{equation}
Combining (\ref{BD}) with the Pollacek--Khintchine formula (\ref{KP}),
we obtain
%
\begin{equation}
\label{BD1} \lim_{u\to\infty} \frac{\Vbar(u-y)}{\Vbar(u)}=e^{\ga y},
\end{equation}
and hence the first result in (\ref{vcon}) holds as claimed. The second
result in (\ref{vcon}) is a consequence of (4.15) in \cite{GMS}, for
example. Since $\Vbar(\ln x)$ is regularly varying, note that the
convergence in (\ref{BD1}) is uniform on compact subsets of $\R$; see
Theorem~1.2.1 of \cite{BGT}.

Let
%
\begin{equation}
\label{VV*} V^*(x):=\int_0^\infty
P^*(H_s\le x)\,ds=\int_{y\le x} e^{\ga y}V(dy),
\end{equation}
see \cite{BD94} or Section~7.2 of \cite{kypbook}. Then $V^*$ is a
renewal function, and so by the key renewal theorem
%
\begin{equation}
\label{KRT} \int_0^u g(y) V^*(u-dy)\to
\frac{1}{m^*}\int_0^\infty g(y)\,dy,
\end{equation}
if $g\ge0$ is directly Riemann integrable on $[0,\infty)$.
We will make frequent use of the following criterion for direct Riemann
integrability.
\textit{If $g\ge0$ is continuous a.e. and dominated by a bounded},
\textit{nonincreasing integrable function on}
$[0,\infty)$, \textit{then} $g$ \textit{is directly Riemann integrable on $[0,\infty)$.}
See Chapter V.4 of \cite{Ab} for information about the key renewal
theorem and direct Riemann integrability.

The function to which we would like to apply (\ref{KRT}), namely
$e^{\ga y}h(y)$ where $h$ is given by (\ref{h}), is typically unbounded
at $0$. To overcome this difficulty, we use the following result.

\begin{prop}\label{concon} If (\ref{Cr}) holds, $h\ge0$, $e^{\ga
y}h(y)1_{[\ve,\infty)}(y)$ is directly Riemann integrable for every
$\ve
>0$, and
%
\begin{equation}
\label{KRTve} \limsup_{\ve\to0}\limsup_{u\to\infty}
\int_{[0,\ve)}h(y)\frac
{V(u-dy)}{\Vbar(u)}=0,
\end{equation}
then (\ref{Fcon}) holds.
\end{prop}

\begin{pf}  By (\ref{KP}) and (\ref{VV*}),
\[
\frac{V(u-dy)}{q\Vbar(u)}= \frac{e^{\ga y}V^*(u-dy)}{e^{\ga u}P(\tau(u)
<\infty)}.
\]
Hence, by (\ref{BD}), and (\ref{KRT}) applied to $e^{\ga
y}h(y)1_{[\ve
,\infty)}(y)$,
\begin{eqnarray*}
\int_{\ve}^u h(y) \frac{V(u-dy)}{q\Vbar(u)} = \int
_{\ve}^u e^{\ga
y}h(y) \frac{V^*(u-dy)}{e^{\ga u}P(\tau(u)<\infty)} &
\to&\int_{\ve}^\infty e^{\ga y} h(y)
\frac{\ga}{q}\,dy
\end{eqnarray*}
as $u\to\infty$. Combined with (\ref{KRTve}) and monotone convergence,
this proves (\ref{conh}), which in turn is equivalent to (\ref{Fcon}).
\end{pf}

The next result gives conditions on $F$ which ensure that $h$, defined
by (\ref{h}), satisfies the hypotheses of Proposition \ref{concon}.

\begin{prop}\label{hDRI} Assume $F\ge0$ satisfies the hypotheses of
Proposition \ref{hcont}, and further that $EX_1e^{\ga X_1}<\infty$.
Then $h$ satisfies the hypotheses of Proposition~\ref{concon}.
\end{prop}

\begin{pf}  For any $y\ge0$,
%
\begin{eqnarray}
\label{hPi} e^{\ga y}h(y)&=& \int_{{\ceb}}e^{\ga y}F(y,
\gep) \nbar\bigl(d\gep, \gep (\gz )> y\bigr)\nonumber
\\
&=& \int_{\ceb} I\bigl(\gep(\gz)>y\bigr)e^{-\ga(\gep(\gz)- y)}F(y,
\gep) e^{\ga
\gep(\gz)}\nbar(d\gep)
\\
&\le& C\int_{(y,\infty)} e^{\ga x}\Pi_H(dx),\nonumber
\end{eqnarray}
by (\ref{nV1}). Further,
\[
\int_{y\ge0}\int_{x> y} e^{\ga x}
\Pi_H(dx) \,dy= \int_{x\ge0} x e^{\ga x}
\Pi_H(dx)<\infty,
\]
by an argument analogous to Proposition 7.1 of \cite{G}. Thus, $e^{\ga
y}h(y)$ is dominated by a nonincreasing integrable function on
$[0,\infty)$, and hence, for each $\ve>0$, $e^{\ga y}h(y)1_{[\ve
,\infty
)}(y)$ is dominated by a bounded nonincreasing integrable function on
$[0,\infty)$. Additionally, by Proposition \ref{hcont}, $h$ is
continuous a.e. with respect to Lebesgue measure. Consequently, $e^{\ga
y}h(y)1_{[\ve,\infty)}(y)$ is directly Riemann integrable for every
$\ve>0$.

Next, since the convergence is uniform on compact in (\ref{BD1}), for
any $x\ge0$,
%
\begin{equation}
\label{VV1} \int_{[0,x)}\frac{V(u-dy)}{\Vbar(u)}\to
e^{\ga x}-1
\end{equation}
as $u\to\infty$.
Thus, by (\ref{hPi}) and (\ref{VV1}), if $\ve<1$ and $u$ is
sufficiently large
\begin{eqnarray*}
&& \int_{[0,\ve)}h(y) \frac{V(u-dy)}{\Vbar(u)}
\\
&&\qquad \le C\int_{[0,\ve)}\frac{V(u-dy)}{\Vbar(u)}\int_{(y,\infty
)}e^{\ga
x}
\Pi_H(dx)
\\
&&\qquad =C\int_{[0,\ve)}e^{\ga x}\Pi_H(dx)\int
_{[0,x)}\frac
{V(u-dy)}{\Vbar
(u)}
\\
&&\quad\qquad{}+ C\int_{[\ve,\infty)}e^{\ga x}\Pi
_H(dx)\int_{[0,\ve)}\frac{V(u-dy)}{\Vbar(u)}
\\
&&\qquad \le C_1 \int_{[0,\ve)}xe^{\ga x}
\Pi_H(dx) + C_1\ve\int_{[\ve
,\infty
)}e^{\ga x}
\Pi_H(dx)
\\
&&\qquad \le C_1e^{\ga} \int_{[0,\ve)}x
\Pi_H(dx)+ C_1e^{\ga} \ve\pibar
_H(\ve) + C_1\ve\int_{[1,\infty)}e^{\ga x}
\Pi_H(dx).
\end{eqnarray*}
Now $\int_{[0,1)}x\Pi_H(dx)<\infty$, since $H$ is a subordinator, thus
$\int_{[0,\ve)}x\Pi_H(dx)\to0$ and $\ve\pibar_H(\ve)\to0$ as
$\ve\to
0$. Combined with $\int_{x\ge1} e^{\ga x}\Pi_H(dx)<\infty$, this
shows that the final expression approaches $0$ as $\ve\to0$.
\end{pf}

As a consequence, we have
the following.

\begin{pf*}{Proof  of Theorem \ref{CrThm}}  This follows immediately from Propositions \ref{concon} and \ref{hDRI}.
\end{pf*}

The convergence in Theorem \ref{CrThm} may alternatively be expressed
in terms of the overshoot $X_{\tau(u)}-u$ rather than the undershoot
of the
maximum $u-\Xbar_{\tau(u)-}$.

\begin{teo}\label{GCr}
Assume $G\dvtx  {\ceb}\times[0,\infty)\to[0,\infty)$ is product
measurable, $e^{-\ga x}G(\gep,x)$ is bounded in $(\gep,x)$ and
$G(\gep
,\cdot )$ is continuous for a.e. $\gep$ w.r.t. $\nbar$. Then under
(\ref{Cr}),
%
\begin{eqnarray}
\label{GCr1} && E^{(u)}G(Y_u, X_{\tau(u)}-u)\nonumber
\\
&&\qquad \to\int_{[0,\infty) }\frac{\ga}q e^{\ga y}\,dy \int
_{{\ceb}\times
(0,\infty)} G(\gep,x) \nbar\bigl(d\gep, \gep(\gz)\in y+dx\bigr)
\\
&&\quad\qquad{} + \rmd _H\frac{\ga
}q G({\mathbf0},0).\nonumber
\end{eqnarray}
\end{teo}

\begin{pf}  Let $F(y,\gep)=G(\gep, \gep(\gz)-y)I(\gep(\gz)\ge y)$. Then $F$
satisfies the conditions of Theorem \ref{CrThm}, and
\[
G(Y_u, X_{\tau(u)}-u)= F(u-\Xbar_{\tau(u)-},
Y_u)
\]
on $\{\tau(u)<\infty\}$.
Consequently, (\ref{MR}) yields
\begin{eqnarray*}
&& E^{(u)} G(Y_u, X_{\tau(u)}-u)
\\
&&\qquad \to\int_{[0,\infty) }\frac{\ga}q e^{\ga y}\,dy \int
_{{\ceb
}}G\bigl(\gep,\gep (\gz)-y\bigr) \nbar\bigl(d\gep, \gep(
\gz)>y\bigr)+\rmd_H\frac{\ga}q G({\mathbf 0},0)
\\
&&\qquad =\int_{[0,\infty) }\frac{\ga}q e^{\ga y}\,dy \int
_{{\ceb}\times
(0,\infty)} G(\gep, x) \nbar\bigl(d\gep, \gep(\gz)\in y+dx\bigr)+
\rmd _H\frac{\ga
}q G({\mathbf0},0),
\end{eqnarray*}
completing the proof.
\end{pf}


\section{Proofs under the convolution equivalent condition}\label{s6}

We begin with the definition of the class ${\mathcal S}^{(\alpha)}$.
As mentioned previously,
we will restrict ourselves to the nonlattice case, with
the understanding that the alternative can be handled by obvious
modifications.
A distribution $F$
on $ [0, \infty)$ with tail $\overline{F}=1-F$ belongs to
the class ${\mathcal S}^{(\alpha)}$, $\alpha> 0$,
if $\overline{F}(u)>0$ for all $u>0$,
%
\begin{equation}
\label{Salph} \lim_{u \to\infty} \frac{\overline{F}(u+x)}{\overline{F}(u)} =
e^{-\alpha x}\qquad\mbox{for } x \in(-\infty, \infty),
\end{equation}
and
%
\begin{equation}
\label{S2} \lim_{u \to\infty} \frac{\overline{F^{2*}}(u)}{\overline{F}(u)}\qquad\mbox{exists and
is finite},
\end{equation}
where $F^{2*}=F * F$.
Distributions in ${\mathcal S}^{(\alpha)}$ are called
\textit{convolution equivalent} with index $\alpha$.
When $F\in{\mathcal S}^{(\alpha)}$, the limit in (\ref{S2})
must be of the form $ 2\delta_{\alpha}^F$,
where $\delta_{\alpha}^F:= \int_{[0, \infty)} e^{\alpha
x}F(d x)$ is finite.
Much is known about the properties of such distributions; see, for
example, \cite{C,EG,kl,Pakes,Pakes2}
and \cite{Wat}. In particular, the class is closed under tail equivalence,
that is, if $F\in{\mathcal S}^{(\alpha)}$ and $G$ is a
distribution function for which
\[
\lim_{u\to\infty}\frac{\overline{G}(u)}{\overline{F}(u)}=c \qquad\mbox {for some } c
\in(0,\infty),
\]
then $G\in{\mathcal S}^{(\alpha)}$.

The\vspace*{2pt} convolution equivalent model (\ref{Sa}) was introduced by Kl\"
{u}ppelberg, Kyprianou and Maller \cite{kkm}.\hskip.2pt\footnote{In \cite{kkm},
(\ref{Sa}) is stated in terms of $\Pi_X^+(\cdot\cap[1,\infty))/\Pi
_X^+([1,\infty))\in\mathcal{S}^{ (\alpha )}$. This is
equivalent to $X_1^+\in\mathcal{S}^{ (\alpha )}$ by Watanabe
\cite{Wat}.\label{ftnP}} As noted earlier,
when (\ref{Sa}) holds, $Ee^{\gt X_1}=\infty$ for all $\gt>\ga$, so
(\ref{Cr}) must fail.
Nevertheless, (\ref{vcon}) continues to hold under (\ref{Sa}). This is
because by (\ref{Vinf}), $F(u)=qV(u)$ is a distribution function, and
combining several results in \cite{kkm} (see (4) of \cite{DK}),
together with closure of ${\mathcal S}^{(\alpha)}$ under tail
equivalence, it follows that $F\in{\mathcal S}^{(\alpha)}$. Hence, the
first condition in (\ref{vcon}) follows from (\ref{Salph}). The second
condition, which corresponds to asymptotic creeping, again follows from
results in \cite{kkm} and can also be found in \cite{DK}.

We begin with a general result about convolution equivalent distributions.

\begin{lemma}\label{conhF}
If $F\in{\mathcal S}^{(\alpha)}$, and $g\ge0$ is continuous
a.e. (Lebesgue) with $g(y)/\Fbar(y)\to L$ as $y\to\infty$, then
\[
\label{hFlim} \int_{0\le y\le u} g(y) \frac{F(u-dy)}{\Fbar(u)}\to\int
_{0}^\infty g(y) \ga e^{\ga y}\,dy+L\int
_{0}^\infty e^{\ga y} F(dy) \qquad\mbox{as }
u\to\infty.
\]
\end{lemma}

\begin{pf}  Fix $K\in(0,\infty)$ and write
%
\begin{eqnarray}
\label{hFlim1}
&& \int_{0\le y\le u} g(y) \frac{F(u-dy)}{\Fbar(u)}\nonumber
\\
&&\qquad =  \biggl(\int_{0\le
y\le K}+\int_{K< y<u-K}+\int
_{u-K\le y\le u} \biggr) g(y)
\frac{F(u-dy)}{\Fbar(u)}
\\
&&\qquad =\mathrm{I}+\mathrm{II}+\mathrm{III}.\nonumber
\end{eqnarray}
By vague convergence,
\[
\mathrm{I}\to\int_{0}^K g(y) \ga
e^{\ga y}\,dy.
\]
Next,
\[
\mathrm{III}= \int_{0\le y\le K} g(u-y) \frac{F(dy)}{\Fbar(u)}= \int
_{0\le
y\le K} \frac{g(u-y)}{\Fbar(u-y)} \frac{\Fbar(u-y)}{\Fbar(u)}F(dy).
\]
For large $u$, the integrand is bounded by $2Le^{\ga K}$ and converges
to $Le^{\ga y}$, thus by bounded convergence,
\[
\mathrm{III}\to L\int_{0}^K e^{\ga y}
F(dy).
\]
Finally,
\[
\limsup_{K\to\infty}\limsup_{u\to\infty} \mathrm{II} \le
\limsup_{K\to
\infty}\limsup_{u\to\infty}\sup
_{y\ge K}\frac{g(y)}{\Fbar
(y)}\int_{K<
y<u-K} {
\Fbar(y)}\frac{F(u-dy)}{\Fbar(u)}=0,
\]
by Lemma 7.1 of \cite{kkm}. Thus, the result follows by letting $u\to
\infty$ and then $K\to\infty$ in (\ref{hFlim1}).
\end{pf}

We now turn to conditions under which (\ref{Fcon}) holds in terms of
$h$ given by (\ref{h}).

\begin{prop}\label{conconSa}
If (\ref{Sa}) holds, and $h\ge0$ is continuous a.e. (Lebesgue) with
$h(y)/\Vbar(y)\to0$ as $y\to\infty$, then (\ref{Fcon}) holds. More
generally, assume
$h(y)/\Vbar(y)\to L$, then an extra term needs to be added to the RHS
of (\ref{Fcon}), namely ${L}/{q\gk(0,-\ga)}$.
\end{prop}

\begin{pf}  As noted above, $qV(u)$ is a distribution function in ${\mathcal
S}^{(\alpha)}$. Thus, by Lemma \ref{conhF},
\[
\int_{0\le y\le u} h(y) \frac{V(u-dy)}{\Vbar(u)}\to\int
_{0}^\infty h(y) \ga e^{\ga y}\,dy+L\int
_{0}^\infty e^{\ga y} V(dy).
\]
Dividing through by $q$ and using (\ref{Vkap}) and (\ref{kpos}) gives
\[
\int_{0\le y\le u} h(y) \frac{V(u-dy)}{q\Vbar(u)}\to\int
_{0}^\infty h(y) \frac{\ga e^{\ga y}}{q}\,dy+
\frac{L}{q\gk(0,-\ga)}.
\]
With $L=0$, this is (\ref{conh}) which is equivalent to (\ref{Fcon}).
\end{pf}

The next result gives conditions on $F$ in (\ref{h}) which ensures
convergence of $h(y)/\Vbar(y)$ as $y\to\infty$.

\begin{prop}\label{Saprop}
If (\ref{Sa}) holds and
\[
\label{FSa1} \bigl\llvert F(y,\gep)-L\bigr\rrvert I\bigl(\gep(\gz)>y\bigr) \to0
\qquad\mbox{uniformly in } \gep \in {\ceb} \mbox{ as } y\to\infty,
\]
then $h(y)/\Vbar(y)\to L\gk^2(0,-\ga)$.
\end{prop}

\begin{pf}  By (\ref{h}),
\[
h(y)\sim L \nbar\bigl(\gep(\gz)>y\bigr)=L\pibar_H(y)\sim L
\gk^2(0,-\ga )\Vbar(y)
\]
by (\ref{KP}) together with (4.4) of \cite{kkm}.
\end{pf}

As a consequence, we have
the following.

\begin{pf*}{Proof  of Theorem \ref{SaThm}}  This follows immediately from Propositions \ref{conconSa} and \ref{Saprop}.
\end{pf*}

\begin{remark}\label{RemF}
Another condition under which (\ref{FSa}) from Theorem \ref{SaThm}
holds, other than when $F$ has compact support, is if $F(y,\gep
)=\tilde
{F}(y,\gep)I(\gep(\gz)\le K)$ for some function $\tilde{F}$ and some
$K\ge0$. In particular if $\tilde{F}\ge0 $ satisfies the hypotheses of
Proposition \ref{hcont}, then $F$ satisfies all the hypotheses of
Theorem \ref{SaThm}.
\end{remark}

The convergence in (\ref{Qconv}) cannot be improved to
$\stackrel{\mathrm{w}}{\longrightarrow}$ since from (\ref{Qinf})
the total mass of $Q^{(\infty)}$ is given by
%
\begin{eqnarray}
\label{Qmass} \bigl\llvert Q^{(\infty)}\bigr\rrvert &=& \frac{\ga}q
\int_{[0,\infty)\times{\ceb}
}e^{\ga y}\nbar \bigl(d\gep, \gep(\gz)>y
\bigr)\,dy+ \rmd_H\frac{\ga}q\nonumber
\\
&=& \frac{\ga}q \int_{[0,\infty)}e^{\ga y}
\pibar_H(y)\,dy + \rmd _H\frac{\ga
}q\nonumber\\[-8pt]\\[-8pt]\nonumber
&=& \frac{1}q \int_{[0,\infty)}\bigl(e^{\ga y}-1
\bigr)\Pi_H(dy) + \rmd_H\frac
{\ga
}q
\\
&=& 1-\frac{\gk(0,-\ga)}q.\nonumber
\end{eqnarray}
Under (\ref{Cr}), $\gk(0,-\ga)=0$ so $\llvert  Q^{(\infty)}\rrvert  =1$, but under
(\ref{Sa}),
$\gk(0,-\ga)>0$ and so $\llvert  Q^{(\infty)}\rrvert  <1$.

As with Theorem \ref{CrThm}, the convergence in Theorem \ref{SaThm} may
alternatively be expressed in terms of the overshoot $X_{\tau(u)}-u$ rather
than the undershoot $u-\Xbar_{\tau(u)-}$.

\begin{teo}\label{GSa}
Assume $G\dvtx  {\ceb}\times[0,\infty)\to[0,\infty)$ is product
measurable, $e^{-\ga x}G(\gep,x)$ is bounded in $(\gep,x)$ and
$G(\gep
,\cdot )$ is continuous for a.e. $\gep$ w.r.t. $\nbar$. Further
assume that
\[
\label{GSa1} G\bigl(\gep,\gep(\gz)-y\bigr)I\bigl(\gep(\gz)>y\bigr)\to0 \qquad
\mbox{uniformly in } \gep \in{\ceb} \mbox{ as } y\to\infty.
\]
Then under (\ref{Sa}),
%
\begin{eqnarray}
\label{GCE}
&& E^{(u)} G(Y_u, X_{\tau(u)}-u)\nonumber
\\
&&\qquad \to\int_{[0,\infty) }\frac{\ga}q e^{\ga y}\,dy \int
_{{\ceb}\times
(0,\infty)} G(\gep,x) \nbar\bigl(d\gep, \gep(\gz)\in y+dx\bigr)
\\
&&\quad\qquad{} +\rmd _H\frac{\ga
}q G({\mathbf0},0).\nonumber
\end{eqnarray}
\end{teo}


\section{The irregular case}\label{s7}

We briefly consider the special case of Theorems \ref{GCr} and \ref
{GSa} where 0 is irregular for $(0,\infty)$ for $X$.
In addition to covering the natural L\'evy process version of
Asmussen's random walk result (\ref{As}), that is when $X$ is compound
Poisson, it also includes the widely studied \textit{compound Poisson
model}, which recall includes a negative drift.
We begin by identifying $\nbar$ in terms of the stopped process
$X_{[0,\tau(0)]}$ where
%
\begin{equation}
\label{Xpath} X_{[0,\tau(0)]}(t):=X_{t\wedge\tau(0)},\qquad t\ge0.
\end{equation}

\begin{prop}\label{nbarX}
Assume 0 is irregular for $(0,\infty)$ for $X$, then
%
\begin{equation}
\label{n1d} P\bigl(\tau(0)<\infty\bigr) \nbar(d\gep)=\llvert
\Pi_H\rrvert P(X_{[0,\tau(0)]}\in d\gep).
\end{equation}
\end{prop}

\begin{pf}  By construction, or using the compensation
formula as in Theorem~\ref{ThmQ}, for some constant $c\in(0,\infty)$,
%
\begin{equation}
\label{n1a} \nbar(d\gep)=n(d\gep)=c P(X_{[0,\tau(0)]}\in d\gep).
\end{equation}
Since $P(X_{\tau(0)}=0, \tau(0)<\infty)=0$, this implies
%
\begin{equation}
\label{n1b} c P\bigl(\tau(0)<\infty\bigr)=\nbar\bigl(\gep(\gz)>0,\gz<\infty
\bigr)=\llvert \Pi_H\rrvert
\end{equation}
by (\ref{nV1}).
Combining (\ref{n1a}) and (\ref{n1b}) proves (\ref{n1d}).
\end{pf}

\begin{prop}\label{0irreg}
Assume 0 is irregular for $(0,\infty)$ for $X$ and either, $G$ is as in
Theorem \ref{GCr} and (\ref{Cr}) holds, or $G$ is as in Theorem \ref
{GSa} and (\ref{Sa}) holds, then
%
\begin{eqnarray}
\label{0irr}
&& E^{(u)} G(Y_u, X_{\tau(u)}-u)
\nonumber\\[-8pt]\\[-8pt]\nonumber
&&\qquad \to\frac{\ga\llvert  \Pi_H\rrvert  \int_0^\infty e^{\ga y}E\{G(X_{[0,\tau(0)]},
X_{\tau(0)}-y); X_{\tau(0)}>y, \tau(0)<\infty\} \,dy}{q P(\tau
(0)<\infty)}.\hspace*{-15pt}
\end{eqnarray}
\end{prop}

\begin{pf}  Since (\ref{Cr}) or (\ref{Sa}) holds, we have $P(\tau(0)<\infty)>0$.
Thus, by (\ref{n1d}), if $y\ge0, x\ge0$, then
\[
\label{n} \nbar\bigl(d\gep, \gep(\gz)\in y+dx\bigr)=\frac{\llvert  \Pi_H\rrvert  P(X_{[0,\tau
(0)]}\in d\gep
, X_{\tau(0)}\in y+dx, \tau(0)<\infty)}{P(\tau(0)<\infty)}.
\]
Since $H$ is compound Poisson when $0$ is irregular for $(0,\infty)$,
we have $\rmd_H=0$. Consequently
(\ref{GCr1}) or (\ref{GCE}) yields
\[
E^{(u)}G(Y_u, X_{\tau(u)}-u)\to\int
_{[0,\infty) }\frac{\ga}q e^{\ga y}\,dy \int
_{{\ceb}\times(0,\infty)} G(\gep,x) \nbar\bigl(d\gep, \gep(\gz)\in y+dx\bigr),
\]
which, by (\ref{n}), is equivalent to (\ref{0irr}).
\end{pf}

Proposition \ref{0irreg} thus provides a natural L\'evy process version
of (\ref{As}) under (\ref{Sa}) as well as under (\ref{Cr}).
We conclude this section by confirming that the constants outside the
integrals in (\ref{As}) and (\ref{0irr}) are in agreement when (\ref
{Cr}) holds. To be precise, by (\ref{BD}), the natural L\'evy process
form of the constant in (\ref{As}), when (\ref{Cr}) holds, is
\[
\frac{\ga m^*}{q E(X_{\tau(0)} e^{\ga X_{\tau(0)}};\tau(0)<\infty)}.
\]
To see that this agrees with the constant in (\ref{0irr}), it suffices
to prove the following.

\begin{lemma}\label{irr}
If 0 is irregular for $(0,\infty)$ for $X$ and (\ref{Cr}) holds, then
\[
\label{id} \llvert \Pi_H\rrvert E\bigl(X_{\tau(0)}
e^{\ga X_{\tau(0)}};{\tau(0)}<\infty\bigr)=P\bigl(\tau (0)<\infty
\bigr)E^*H_1.
\]
\end{lemma}

\begin{pf}  By (\ref{nV1}) and (\ref{n1d}),
\[
\label{Xtau0} \llvert \Pi_H\rrvert P\bigl(X_{\tau(0)}\in
dx, {\tau(0)}<\infty\bigr)=P\bigl(\tau(0)<\infty \bigr)\Pi_{{H}}(dx),
\]
and so
\[
\label{P1} \llvert \Pi_H\rrvert E\bigl(X_{\tau(0)}e^{\ga X_{\tau(0)}};{
\tau(0)}<\infty\bigr)=P\bigl(\tau (0)<\infty\bigr) \int_0^\infty
x e^{\ga x} \Pi_H(dx).
\]
Since $\mathrm{d}_H=0$ when $0$ is irregular for $(0,\infty)$, the result now
follows from (\ref{m*}).
\end{pf}


\section{Proofs of Theorems \texorpdfstring{\protect\ref{thmf}}{3.3},
\texorpdfstring{\protect\ref{thmfSa}}{3.4} and related results}\label{s8}

To calculate the limits in Theorems \ref{thmf} and \ref{thmfSa}, we
consider a particular form for $F$ in Theorems \ref{CrThm} and~\ref{SaThm}. Let $f\dvtx [0,\infty)^4\to[0,\infty)$ be a Borel function, and set
%
\begin{equation}
\label{Fov} F(y,\gep)=f\bigl(y, \gep(\gz)-y,y- \gep(\gz-),\gz\bigr)I\bigl(
\gep(\gz)\ge y\bigr).
\end{equation}
Then
\[
F(u-\Xbar_{\tau(u)-}, Y_u)=f\bigl(u-\Xbar_{\tau(u)-},
X_{\tau(u)}-u, u-X_{\tau(u)-}, \tau(u)- G_{\tau(u)-}\bigr)
\]
on $\{\tau(u)<\infty\}$.
To calculate the limit in this case, we need the following.

\begin{lemma}\label{LGf} If $F$ is of the form (\ref{Fov}) then for
every $y\ge0$,
%
\begin{eqnarray}
\label{Gf}
\qquad && \int_{{\ceb}}F(y,\gep) \nbar\bigl(d\gep, \gep(
\gz)>y\bigr)
\nonumber\\[-8pt]\\[-8pt]\nonumber
&&\qquad =\int_{x>0}\int_{v\ge0}\int
_{t\ge0} f(y,x,v,t)I(v\ge y)
\whV(dt,dv-y)\Pi_X(v+dx).
\end{eqnarray}
If in addition, $Ee^{\ga X_1}<\infty$,
$f$ is jointly continuous in the first three variables
and
$e^{-\ga x}f(y,x,v,t)$ is bounded, then $F$ satisfies the hypotheses
of Proposition \ref{hcont}.
\end{lemma}

\begin{pf}  Using Proposition \ref{nmarg} in the third equality, we have
\begin{eqnarray*}
&& \int_{{\ceb}} F(y,\gep) \nbar\bigl(d\gep, \gep(\gz)>y\bigr)
\\
&&\qquad =\int_{{\ceb}}f\bigl(y, \gep(\gz)-y,y- \gep(\gz-),\gz\bigr)
\nbar\bigl(d\gep, \gep(\gz )>y\bigr)
\\
&&\qquad =\int_{x>0}\int_{z\ge0}\int
_{t\ge0} f(y,x,y+z,t) \nbar\bigl(\gep(\gz )\in y+ dx,\gep(\gz-)
\in-dz, \gz\in dt \bigr)
\\
&&\qquad =\int_{x>0}\int_{z\ge0}\int
_{t\ge0} f(y,x,y+z,t)\whV(dt,dz)\Pi _X(z+y+dx)
\\
&&\qquad =\int_{x>0}\int_{v\ge0}\int
_{t\ge0} f(y,x,v,t)I(v\ge y) \whV (dt,dv-y)\Pi_X(v+dx),
\end{eqnarray*}
which proves (\ref{Gf}).

For the second statement, we only need to check $\nbar(B_y^c)=0$ for
a.e. $y$. But
$B_y^c\subset\{\gep\dvtx \gep(\gz)=y\}$, and so $\nbar(B_y^c)\le\Pi
_H(\{y\}
)=0$ except for at most countably many $y$.
\end{pf}

\begin{remark}\label{rem1} Lemma \ref{LGf} remains true if $f$ is
replaced by $\phi(y)f(y,x,v,t)$ where $\phi$ is bounded and
continuous a.e.,
since in that case $\nbar(B_y^c)=0$ except when $y$ is a point of
discontinuity of $\phi$ or $\Pi_H(\{y\})=0$.
\end{remark}

For reference below we note that if $e^{-\ga x}f(y,x,v,t)$ is bounded then
%
\begin{eqnarray}
\label{x0}
&& \int_{x=0}\int_{y\ge0}\int
_{v\ge0}\int_{t\ge0} f(y,x,v,t)e^{\ga y}
\,dy\,I(v\ge y)
\nonumber\\[-8pt]\\[-8pt]\nonumber
&&\hspace*{86pt}{}\times  \whV(dt,dv-y)\Pi_X(v+dx)=0,
\end{eqnarray}
since
\begin{eqnarray*}
&& \int_{x=0}\int_{y\ge0}\int
_{v\ge0}\int_{t\ge0} e^{\ga
y}\,dy\, I(v\ge
y) \whV(dt,dv-y)\Pi_X(v+dx)
\\
&&\qquad =\int_{y\ge0}\int_{v\ge0}e^{\ga y}\,dy\,
I(v\ge y) \whV(dv-y)\Pi_X\bigl(\{ v\} \bigr)
\\
&&\qquad =\int_{v\ge0} \whV(dv)\int_{y\ge0}e^{\ga y}
\Pi_X\bigl(\{v+y\}\bigr) \,dy =0.
\end{eqnarray*}

We first consider limiting results for $F$ of the form (\ref{Fov}) in
the Cram\'er--Lundberg setting, beginning with Theorem \ref{thmf}.

\begin{pf*}{Proof  of Theorem \ref{thmf}}   Define $F$ by (\ref
{Fov}). Then by Lemma \ref{LGf}, $F$ satisfies the hypotheses of Theorem
\ref{CrThm}, and hence the result follows from (\ref{MR}), (\ref{Gf})
and~(\ref{x0}).
\end{pf*}

Marginal convergence in each of the first three variables in (\ref
{jtcon}) was shown in \cite{GMS}. Equation (\ref{jtcon}) exhibits the
stronger joint convergence and includes the additional time variable $
\tau(u)- G_{\tau(u)-}$. Note also that in the time variable, there is no
restriction on $f$ beyond bounded, and hence the convergence is
stronger than weak convergence in this variable.

As an illustration of (\ref{thmf1}) we obtain, for any $\gl\le0,
\eta
\le\ga, \rho\le0$ and $\gd\ge0$,
%
\begin{eqnarray}
\label{GSga}
&&E^{(u)}e^{\lambda(u-\Xbar_{\tau(u)-})+\eta(X_{\tau(u)}-u)+\rho
(u-X_{\tau(u)-})-\delta
(\tau(u)- G_{\tau(u)-}) } \nonumber
\\
&&\qquad \to
\int_{x\ge0}\int_{y\ge0}\int
_{v\ge0}\int_{t\ge0} e^{\lambda
y+\eta x+\rho v-\delta t}
\frac{\ga}q e^{\ga y} \,dy\, I(v\ge y)
\\
&&\hspace*{89pt}\qquad\quad{}\times \whV(dt, dv-y)
\Pi_X(v+dx) + \rmd_H\frac{\ga}q.\nonumber
\end{eqnarray}
This gives the future value, at time $G_{\tau(u)-}$, of a Gerber--Shiu
expected discounted penalty\vspace*{1pt} function (EDPF) as $u\to\infty$. The
present value is zero since $\tau(u)\to\infty$ in $P^{(u)}$
probability as
$u\to\infty$.
The limit can be simplified if $\rho=0$. From (\ref{DKV})
and~(\ref{x0}), we obtain
%
\begin{eqnarray}
\label{EDPF} &&E^{(u)}e^{\lambda(u-\Xbar_{\tau(u)-})+\eta(X_{\tau(u)}-u)-\delta
(\tau(u)- G_{\tau(u)-})}\nonumber
\\
&&\qquad  \to\int_{x> 0}\int_{y\ge0}\int
_{t\ge0} e^{\lambda y+\eta
x-\delta t}\frac{\ga}q e^{\ga y} \,dy
\int_{v\ge y} \whV(dt, dv-y)\Pi_X(v+dx)\nonumber
\\
&&\quad\qquad{}  + \rmd_H\frac{\ga}q\nonumber
\\
&&\qquad  =\int_{x> 0}\int_{y\ge0}\int
_{t\ge0} e^{\lambda y+\eta
x-\delta t}\frac{\ga}q e^{\ga y} \,dy
\Pi_{{L}^{-1},{H}}(dt,y+dx) + \rmd_H\frac
{\ga}q
\\
&&\qquad  =\int_{t\ge0}\int_{y\ge0}\int
_{x> y} e^{\lambda y+\eta
(x-y)-\delta t}\frac{\ga}q e^{\ga y} \,dy
\Pi_{{L}^{-1},{H}}(dt,dx)+ \rmd _H\frac{\ga}q\nonumber
\\
&&\qquad  =\frac{\ga}{q(\eta-\gl-\ga)}\int_{t\ge0}\int_{x>0}
\bigl(e^{\eta
x}-e^{(\gl+\ga) x} \bigr)e^{-\gd t}
\Pi_{{L}^{-1},{H}}(dt,dx)+ \rmd _H\frac{\ga}q\nonumber
\\
&&\qquad  =\frac{\ga(\gk(\gd,-(\gl+\ga))-\gk(\gd,-\eta))}{q(\eta-\gl
-\ga)}.\nonumber
\end{eqnarray}
\smallskip

Under (\ref{Cr}), it is possible that $Ee^{\gt X_1}=\infty$ for all
$\gt
>\ga$, but it is often the case that $Ee^{\gt X_1}<\infty$ for some
$\gt
>\ga$. The next result extends Theorem \ref{thmf} to include this
possibility, and also provides more information when restricted to the
former setting. This is done by taking advantage of the special
form of $F$ in (\ref{Fov}), whereas Theorem \ref{thmf} was derived from
the general convergence result in Theorem~\ref{CrThm}. It is
interesting to note how the exponential moments may be spread out over
the undershoot variables.
The EDPF results in (\ref{GSga}) and (\ref{EDPF}) also have obvious
extensions to this setting.

\begin{teo}\label{thmft} Assume (\ref{Cr}) holds and $f\dvtx [0,\infty
)^4\to
[0,\infty)$ is a Borel function which is jointly continuous in the
first three variables. Assume $\gt\ge\ga$ and one of the following
three conditions holds:

\begin{longlist}[(iii)]
\item[(i)] $Ee^{\gt X_1}<\infty$, $\rho<\gt$ and $\lambda+ \rho< \gt-\ga
$;

\item[(ii)] $EX_1e^{\gt X_1}<\infty$, $\rho\le\gt$ and $\lambda+ \rho\le
\gt
-\ga$, with at least one of these inequalities being strict;

\item[(iii)] $EX^2_1e^{\gt X_1}<\infty$, $\rho\le\gt$ and $\lambda+ \rho
\le
\gt-\ga$.
\end{longlist}

If $e^{-\lambda y-\gt x-\rho v}f(y,x,v,t)$ is bounded, then
%
\begin{eqnarray}
\label{ft1} && E^{(u)}f\bigl( u-\Xbar_{\tau(u)-},X_{\tau(u)}-u,
u-X_{\tau(u)-}, \tau (u)- G_{\tau(u)-}\bigr)\nonumber
\\
&&\qquad  \to\int_{x\ge0}\int_{y\ge0}\int
_{v\ge0}\int_{t\ge0} f(y,x, v,t)
\frac{\ga}q e^{\ga y} \,dy\, I(v\ge y)
\nonumber\\[-8pt]\\[-8pt]\nonumber
&&\hspace*{122pt}{}\times {} \whV(dt,dv-y)
\Pi_X(v+dx)
\\
&&\quad\qquad{}+\rmd_H\frac{\ga}q f(0,0,0,0). \nonumber
\end{eqnarray}
\end{teo}

\begin{pf}   Define $F$ by (\ref{Fov}) and then $h$ by
(\ref{h}). We will show that $h$ satisfies the hypotheses of
Proposition \ref{concon}. Let $\tf(y,x,v,t)=e^{-\lambda y-\gt x-\rho
v}f(y,x,v,t)$. Then $\tf$ is bounded, jointly continuous in the first
three variables, and by (\ref{Gf}), for every $y\ge0$,
%
\begin{eqnarray}
\label{Gf1} h(y)&=&\int_{{\ceb}}F(y,e) \nbar\bigl(de, e(
\gz)>y\bigr)\nonumber
\\
&=&\int_{x>0}\int_{v\ge0}\int
_{t\ge0} \tf(y,x,v,t)I(v\ge y)e^{\lambda
y+\gt x+\rho v}\nonumber
\\
&&\hspace*{63pt}{}\times  \whV(dt,dv-y)
\Pi_X(v+dx)
\\
&=&e^{(\lambda+\rho-\gt)y}\int_{x>0}\int_{v\ge0}
\int_{t\ge0} \tf (y,x-v-y,v+y,t)I(x> v+y)\nonumber
\\
&&\hspace*{108pt}{} \times e^{\gt
x+(\rho-\gt)v}
\whV(dt,dv)\Pi_X(dx)\nonumber
\end{eqnarray}
after a change of variables. Let $g_y(x,v,t)= \tf(y,x-v-y,v+y,t)I(x>
v+y)e^{\gt x+(\rho-\gt)v}$. Then clearly $g_z(x,v,t)\to g_y(x,v,t)$
as $z\downarrow y$, for every $y\ge0$.
Now fix $y>0$ and let $\llvert  z-y\rrvert  <y/2$. Then for some constant $C$
independent of $z,x,v$ and $t$,
\[
\bigl\llvert g_z(x,v,t)\bigr\rrvert \le CI(x> v+y/2)e^{\gt x+(\rho-\gt)v}
\]
and
%
\begin{eqnarray}
\label{LG}
&&\int_{x>0}\int_{v\ge0}\int
_{t\ge0} I(x> v+y/2) e^{\gt x+(\rho-\gt
)v} \whV(dt,dv)
\Pi_X(dx)
\nonumber\\[-8pt]\\[-8pt]\nonumber
&&\qquad \le\int_{x> y/2}e^{\gt x}\Pi_X(dx)\int
_{0\le v<x-y/2} e^{(\rho
-\gt)v} \whV(dv).
\end{eqnarray}
If we show this last integral is finite, then by dominated convergence,
$h(z)\to h(y)$ as $z\downarrow y$ for every $y>0$, showing that $h$ is
right continuous on $(0,\infty)$, and consequently continuous a.e. The
final expression in (\ref{LG}) is decreasing in $y$, hence to prove
finiteness it suffices to prove the following stronger result, which
will be needed below; for every $\ve>0$,
%
\begin{eqnarray}
\label{Ive} I_{\ve}&:=&\int_\ve^\infty
e^{(\ga+\lambda+\rho-\gt)y}\,dy\int_{x>
y}e^{\gt x}
\Pi_X(dx) \int_{0\le v <x-y}e^{(\rho-\gt)v} \whV (dv)<
\infty.\hspace*{-30pt}
\end{eqnarray}
We will need the following consequence of Proposition 3.1 of Bertoin
\cite{bert}; for every $y>0$ there is a constant $c=c(y)$ such that
%
\begin{equation}
\label{Vbd} \whV(v)\le cv \qquad\mbox{for } v\ge y.
\end{equation}
First assume $\rho<\gt$. Integrating by parts and using (\ref{Vbd})
shows that the last of the three integrals in (\ref{Ive}) is bounded
independently of $x$ and $y$, hence
\begin{eqnarray*}
I_{\ve}&\le& C \int_{x>\ve}e^{\gt x}
\Pi_X(dx) \int_{\ve<y<x} e^{(\ga
+\lambda+\rho-\gt)y}\,dy
\\
&\le& C \cases{\displaystyle \int_{x>\ve}e^{\gt x}
\Pi_X(dx), &\quad if $\ga+\lambda+\rho-\gt <0$,
\vspace*{5pt}\cr
\displaystyle\int
_{x>\ve}xe^{\gt x}\Pi_X(dx), &\quad if $
\ga+\lambda+\rho-\gt=0$.}
\end{eqnarray*}
Thus, $I_{\ve}<\infty$ under each of the assumptions (i), (ii) and (iii) by Theorem 25.3 of Sato \cite{sato}. Now assume
$\rho=\gt$. Then
\[
I_{\ve}=\int_{x>\ve}e^{\gt x}
\Pi_X(dx) \int_{\ve<y<x} e^{(\ga
+\lambda
+\rho-\gt)y}
\whV(x-y)\,dy.
\]
If $\ga+\lambda+\rho-\gt=0$, then we are in case {(iii)} and
\[
I_{\ve}\le\int_{x>\ve}x\whV(x)e^{\gt x}
\Pi_X(dx)\le C \int_{x>\ve
}x^{2}e^{\gt x}
\Pi_X(dx),
\]
which is finite under {(iii)}. Finally, if $\ga+\lambda+\rho-\gt<0$
then we are in case {(ii)} or~{(iii)}. We break $I_{\ve}$ into
two parts $I_{\ve}(1)+I_{\ve}(2)$ where
\begin{eqnarray*}
I_{\ve}(1)&=&\int_{x>\ve}e^{\gt x}
\Pi_X(dx) \int_{\ve\vee(x-1)<y<x} e^{(\ga+\lambda+\rho-\gt)y}\whV(x-y)\,dy
\\
&\le&\whV(1) \int_{x>\ve}e^{\gt x}\Pi_X(dx)
\end{eqnarray*}
and
\begin{eqnarray*}
I_{\ve}(2)&=&\int_{x>\ve}e^{\gt x}
\Pi_X(dx) \int_{\ve<y\le\ve
\vee
(x-1)} e^{(\ga+\lambda+\rho-\gt)y}\whV(x-y)\,dy
\\
&\le& c(1)\int_{x>\ve}e^{\gt x}\Pi_X(dx)
\int_{\ve<y\le\ve\vee(x-1)} e^{(\ga+\lambda+\rho-\gt)y}(x-y)\,dy
\\
&\le& C\int_{x>\ve}xe^{\gt x}\Pi_X(dx).
\end{eqnarray*}
Thus, $I_{\ve}$ is finite in this case also, completing the proof of
(\ref{Ive}).

By (\ref{Gf1}), for every $y\ge0$,
%
\begin{eqnarray}
\label{hk} e^{\ga y}h(y)&\le& C e^{(\ga+\lambda+\rho-\gt)y}\int
_{v\ge0}\int_{x>0} I(x> v+y)e^{\gt x+(\rho-\gt)v}
\whV(dv)\Pi_X(dx)\hspace*{-15pt}
\nonumber\\[-8pt]\\[-8pt]\nonumber
&=:&k(y)
\end{eqnarray}
say. Clearly, $k$ is nonincreasing on $[0,\infty)$, and for every $\ve>0$,
\begin{eqnarray*}
\int_\ve^\infty k(y)\,dy& =& C\int
_\ve^\infty e^{(\ga+\lambda+\rho
-\gt
)y}\,dy\int
_{x> y}e^{\gt x}\Pi_X(dx) \int
_{0\le v <x-y}e^{(\rho-\gt)v} \whV(dv)
\\
&<&\infty,
\end{eqnarray*}
by (\ref{Ive}) under {(i)}, {(ii)} or {(iii)}.
Hence, in each case $e^{\ga y}h(y)1_{[\vep,\infty)}(y)$ is directly
Riemann integrable for every $\ve>0$.

Finally, from (\ref{hk}), for $\ve\in(0,1)$,
\begin{eqnarray*}
&& \int_{[0,\ve)} h(y)\frac{V(u-dy)}{\Vbar(u)}
\\
&&\qquad \le  C \int
_{v\ge 0}\int_{x>v} e^{\gt x+(\rho-\gt)v} \whV(dv)
\Pi_X(dx)
\int_{[0,\ve)}I(y< x-v)
\frac
{V(u-dy)}{\Vbar(u)}.
\end{eqnarray*}
By the uniform convergence on compact sets in (\ref{BD1}), it follows that
for large $u$,
\begin{eqnarray*}
&& \int_{[0,\ve)} h(y)\frac{V(u-dy)}{\Vbar(u)}
\\
&&\qquad \le  C_1
\biggl(\int_{v\le
1}+\int_{v> 1} \biggr)\int
_{x>v} e^{\gt x+(\rho-\gt)v}
\bigl[(x-v)\wedge\ve\bigr]\whV(dv)\Pi_X(dx)
\\
&&\qquad =\mathrm{I}+\mathrm{II}.
\end{eqnarray*}
Now, using (\ref{EA}),
\begin{eqnarray*}
\mathrm{I}&\le& C_1 \int_{v\le1}\int_{x>0}
(x\wedge\ve) e^{\gt x}\whV (dv)\Pi _X(v+dx)
\\
&\le& C_1 \int_{x>0} (x\wedge
\ve)e^{\gt x}\Pi_H(dx)\to0
\end{eqnarray*}
as $\ve\to0$ by dominated convergence, since $\int_{x>1} e^{\gt
x}\Pi
_H(dx)<\infty$ by Proposition~7.1 of \cite{G} and $\int_{x\le1}
xe^{\gt x}\Pi_H(dx)<\infty$ because $H$ is a subordinator.
For the second term,
\[
\mathrm{II}\le C_1\ve\int_{x>1} e^{\gt x}
\Pi_X(dx) \int_{0\le v<x}e^{(\rho
-\gt
)v} \whV(dv)
\to0
\]
as $\ve\to0$, since the integral is easily seen to be finite from
(\ref{Ive}) and the renewal theorem.
Thus we may apply Proposition \ref{concon} to $h$, and (\ref{ft1})
follows after observing that the integral over $x=0$ in (\ref{ft1}) is
zero by (\ref{x0}).
\end{pf}

\begin{remark} As the proof shows, conditions under which (\ref{ft1})
holds can be stated in terms of the integral condition (\ref{Ive}) on
the renewal function $\whV$, rather than in terms of conditions (i)--(iii).
Specifically, assume (\ref{Cr}) holds, $Ee^{\gt X_1}<\infty$ for some
$\gt\ge\ga$, and $f\dvtx [0,\infty)^4\to[0,\infty)$ is a Borel function
which is jointly continuous in the first three variables and
$e^{-\lambda y-\gt x-\rho v}f(y,x,\break v,t)$ is bounded where $\lambda+
\rho\le\gt-\ga$ and $\rho\le\gt$. If (\ref{Ive}) holds for every
$\gep>0$, then (\ref{ft1}) holds.
\end{remark}

We now turn to the convolution equivalent setting. In this case, we
need to impose an extra condition on $f$ in (\ref{Fov}).

\begin{prop}\label{hfSa}
Assume $F$ is given by (\ref{Fov}) where
%
\begin{equation}
\label{fun0} \sup_{x>0,t\ge0, v\ge y}f(y,x,v,t)\to0 \qquad\mbox{as } y\to
\infty.
\end{equation}
Then (\ref{FSa}) holds.
\end{prop}

\begin{pf}
From (\ref{Fov}),
\begin{eqnarray*}
\sup_{\gep\in{\ceb}}F(y,\gep)I\bigl(\gep(\gz)>y\bigr)&=&\sup
_{\gep\in
{\ceb}}f\bigl(y, \gep(\gz)-y,y- \gep(\gz-),\gz\bigr)I\bigl(\gep(
\gz)>y\bigr)
\\
&\le&\sup_{x>0,t\ge0, v\ge y}f(y,x,v,t)\to0
\end{eqnarray*}
as $y\to\infty$ by (\ref{fun0}).
\end{pf}

As an immediate consequence, we have the following.

\begin{pf*}{Proof of Theorem \ref{thmfSa}}
Define $F$ by (\ref{Fov}). By
Lemma \ref{LGf} and Proposition \ref{hfSa}, $F$ satisfies the
hypotheses of Theorem \ref{SaThm}. Thus, (\ref{Fcon}) holds which is
equivalent to (\ref{Safcon}) by (\ref{Gf}) and (\ref{x0}).
\end{pf*}

Theorem \ref{thmfSa} imposes the extra condition (\ref{fun0}) on $f$
compared with Theorem~\ref{thmf}. As a typical example, any function
$f$ satisfying the conditions of Theorem \ref{thmf}, when multiplied by
a bounded continuous function with compact support $\phi(y)$, trivially
satisfies the conditions of Theorem \ref{thmfSa}. This manifests itself
in the convergence of (\ref{jtconv}) only being vague convergence
rather than weak convergence. It cannot be improved to the
weak convergence of (\ref{jtcon}) since, as noted earlier in (\ref
{Qmass}), the total mass of the limit in (\ref{jtconv}) is $1-\gk
(0,-\ga)q^{-1}$.

Another example of the effect of condition (\ref{fun0}) is in the
calculation of the EDPF analogous (\ref{GSga}). Using Remark \ref
{rem1}, the continuity assumption on $\phi$ above can be weakened to
continuous a.e. Hence, we may take
$\phi(y)=I(y\le K)$ for some $K\ge0$. Thus,
applying (\ref{Safcon}) to the function $f(y,x,v,t)= e^{\lambda y+\eta
x+\rho v-\delta t}I(y\le K)$ where $K>0, \gl\le0, \eta\le\ga, \rho
\le
0$ and $\gd\ge0$, we obtain
\begin{eqnarray*}
\label{GSga1} &&E^{(u)}\bigl\{e^{\lambda(u-\Xbar_{\tau(u)-})+\eta(X_{\tau(u)}-u)+\rho
(u-X_{\tau(u)
-})-\delta(\tau(u)- G_{\tau(u)-})}; u-\Xbar_{\tau(u)-}
\le K\bigr\}
\\
&&\qquad  \to\int_{0\le y\le K}\int_{x\ge0}\int
_{v\ge0}\int_{t\ge0} e^{\lambda y+\eta x+\rho v-\delta t}
\frac{\ga}q e^{\ga y} \,dy\,I(v\ge y)
\\
&&\hspace*{135pt}{} \times \whV(dt, dv-y)\Pi
_X(v+dx) + \rmd_H\frac{\ga}q.
\end{eqnarray*}
The restriction imposed by $K$ cannot be removed as will be apparent
from Proposition \ref{p6} below.
Finally, we point out that there is no extension of Theorem \ref
{thmfSa} to the setting of Theorem \ref{thmft} since $Ee^{\gt
X_1}=\infty$ for all
$\gt>\ga$.

We next address convergence of the marginals in (\ref{jtcon}) and
(\ref{jtconv}).
As indicated in Section~\ref{s4}, some care is needed under (\ref{Sa})
since, in (\ref{jtconv}), the limit of the marginals is not the
marginal of the limits for the case of the overshoot and $\tau(u)-
G_{\tau(u)-}$.
If $F$ is given by (\ref{Fov}) where $f$ depends only on $x$ and $t$,
then, using (\ref{DKV}), (\ref{Gf}) reduces to
%
\begin{eqnarray}
\label{Ffxt}
&& \int_{{\ceb}}F(y,\gep) \nbar\bigl(d\gep, \gep(
\gz)>y\bigr)
\nonumber\\[-8pt]\\[-8pt]\nonumber
&&\qquad =\int_{x>
0}\int_{t\ge
0} f(x,t)
\Pi_{L^{-1},H}(dt,y+dx),\qquad y\ge0.
\end{eqnarray}
In particular, under (\ref{Cr}), by Theorem \ref{thmf}, for $x\ge0,
t\ge0$,
%
\begin{eqnarray}
\label{ltovt}
&& P^{(u)}\bigl(X_{\tau(u)}-u \in dx,\tau(u)-
G_{\tau(u)-}\in dt\bigr)
\nonumber\\[-8pt]\\[-8pt]\nonumber
&&\qquad \stackrel {\mathrm{w}} {\longrightarrow}
\frac{\ga}q\int_{y\ge0} e^{\ga y}
\Pi_{L^{-1},H}(dt,y+dx) \,dy+ \rmd_H\frac{\ga}q
\delta_{(0,0)}(dx,dt).
\end{eqnarray}
Under (\ref{Sa}), the mass of the limit in (\ref{ltovt}) is less than
one. In this case, an extra term appears in the limit. The distribution
of this additional mass and proof of joint weak convergence under (\ref
{Sa}) is given in the following result.

\begin{prop}\label{p6}
Assume (\ref{Sa}) holds and that $f\dvtx [0,\infty)^2\to[0,\infty)$ is a
Borel function which is continuous in the first variable, and $e^{-\gb
x}f(x,t)$ is bounded for some $\gb<\ga$. Then
%
\begin{eqnarray}
\label{Safconw} &&E^{(u)}  f\bigl(X_{\tau(u)}-u, \tau(u)-
G_{\tau(u)-}\bigr)\nonumber
\\
&&\qquad \to\frac{-\Psi_X(-\rmi\ga) }q \int_{x\ge0}\int_{t\ge0}
f(x,t)\ga e^{-\ga x}\,dx\int_{v\ge0} e^{-\ga v}
\whV(dt,dv)
\nonumber\\[-8pt]\\[-8pt]\nonumber
&&\hspace*{2pt}\quad\qquad{} +\frac{\ga}q\int_{x\ge0}\int
_{t\ge0} f(x,t)\int_{y\ge0} e^{\ga y}
\Pi_{L^{-1},H}(dt,y+dx) \,dy
\\
&&\hspace*{2pt}\quad\qquad{}  + \rmd_H\frac{\ga}q f(0,0).\nonumber
\end{eqnarray}
In particular, we have joint convergence; for ${x\ge0}, {t\ge0}$,
%
\begin{eqnarray}
\label{ltovtS}
&& P^{(u)}\bigl(X_{\tau(u)}-u\in dx, \tau(u)-
G_{\tau(u)-}\in dt\bigr)\nonumber
\\
&&\qquad \stackrel{\mathrm{w}} {\longrightarrow}\frac{ -\Psi_X(-\rmi\ga)
}q \ga
e^{-\ga x}\,dx\int_{v\ge0} e^{-\ga v} \whV(dt,dv)
\\
&&\hspace*{9pt}\quad\qquad{}  +\frac{\ga}q\int_{y\ge0} e^{\ga y}\Pi
_{L^{-1},H}(dt,y+dx) \,dy + \rmd_H\frac{\ga}q
\delta_{(0,0)}(dx,dt).\nonumber
\end{eqnarray}
\end{prop}

\begin{pf}
We will use Proposition \ref{conconSa}. Let
\[
\label{Fov1} F(y,\gep)=f\bigl(\gep(\gz)-y,\gz\bigr)I\bigl(\gep(\gz)\ge y
\bigr).
\]
By Lemma \ref{LGf}, $F$ satisfies the hypothesis of Proposition \ref
{hcont}, hence $h$ is continuous a.e.
Next we evaluate the limit of $h(y)/\pibar_{X}(y)$ as $y\to\infty$. By
(\ref{Gf}), for $y\ge0$,
\begin{eqnarray*}
\label{hot} h(y)&=&\int_{{\ceb}}F(y,\gep) \nbar\bigl(d\gep,
\gep(\gz)>y\bigr)
\\
&=&\int_{x>0}\int_{v\ge0}\int
_{t\ge0} f(x,t)I(v\ge y) \whV (dt,dv-y)\Pi _X(v+dx)
\\
&=&\int_{t\ge0}\int_{v\ge0}\whV(dt,dv)\int
_{x>0} f(x, t)\Pi_X(y+v+dx).
\end{eqnarray*}
Observe that for $v\ge0$, from footnote \ref{ftnP} (page~\pageref{ftnP})
and (\ref{Salph}),
\[
\frac{\Pi_X(y+v+dx)}{\pibar_X(y)}\stackrel{\mathrm {w}} {\longrightarrow}\ga e^{-\ga(v+x)}\,dx
\qquad\mbox {on } [0,\infty) \mbox{ as } y\to\infty.
\]
Further, by Potter's bounds (see, e.g., (4.10) of \cite{GM2}), if
$\gamma\in(\gb,\ga)$ then
\[
\frac{\pibar_X(y+v)}{\pibar_X(y)}\le C e^{-\gamma v} \qquad\mbox{if } v\ge0, y\ge1,
\]
where $C$ depends only on $\gamma$. Thus, for any $y\ge1, v\ge0$ and
$K\ge0$,
%
\begin{eqnarray}\label{fbd}
\int_{x>K} e^{\gb x}\frac{\Pi_X(y+v+dx)}{\pibar_X(y)}
&=&\int_{x>K} \gb e^{\gb x} \frac{\pibar_X(y+v+x)}{\pibar_X(y)}\,dx\nonumber
\\
&&{} +
\frac{e^{\gb K}\pibar
_X(y+v+K)}{\pibar
_X(y)}
\nonumber\\[-8pt]\\[-8pt]\nonumber
&\le& C\int_{x>K} \gb e^{\gb x} e^{-\gamma(v+x)} \,dx+C
e^{\gb
K}e^{-\gamma(v+K)}
\\
&\le& Ce^{-\gamma v-(\gamma-\gb)K}.\nonumber
\end{eqnarray}
Now, for any $v\ge0$ and $K\ge0$ write
%
\begin{eqnarray}
\label{splK}
&& \int_{x>0} f(x, t)\frac{\Pi_X(y+v+dx)}{\pibar_X(y)}\nonumber
\\
&&\qquad =
\biggl(\int_{0<x\le
K}+\int_{x>K} \biggr) f(x,
t)\frac{\Pi_X(y+v+dx)}{\pibar_X(y)}
\\
&&\qquad = \mathrm{I}+\mathrm{II}.\nonumber
\end{eqnarray}
By weak convergence,
\[
I\to\int_{0<x\le K}f(x, t) \ga e^{-\ga(v+x)}\,dx \qquad\mbox{as
} y\to \infty,
\]
and by monotone convergence,
\[
\int_{0<x\le K}f(x, t) \ga e^{-\ga(v+x)}\,dx\to\int
_{x\ge0}f(x, t) \ga e^{-\ga(v+x)}\,dx \qquad\mbox{as } K\to
\infty.
\]
On the other hand, by (\ref{fbd}),
\[
\mathrm{II}\le Ce^{-\gamma v-(\gamma-\gb)K}\qquad\mbox{for } y\ge1.
\]
Thus, letting $y\to\infty$ then $K\to\infty$ in (\ref{splK}) gives
\[
\label{splK1} \int_{x>0} f(x, t)\frac{\Pi_X(y+v+dx)}{\pibar_X(y)}\to\int
_{x\ge
0}f(x, t) \ga e^{-\ga(v+x)}\,dx.
\]
Further, by (\ref{fbd}) with $K=0$, for every $v\ge0$
\[
\int_{x>0} f(x,t)\frac{\Pi_X(y+v+dx)}{\pibar_X(y)} \le Ce^{-\gamma v}.
\]
Hence, by dominated convergence,
\[
\frac{h(y)}{{\pibar_X(y)}}\to\int_{x\ge0}\int_{t\ge0}
f(x,t)\ga e^{-\ga x}\,dx\int_{v\ge0} e^{-\ga v}
\whV(dt,dv).
\]
Since
\[
\lim_{y\to\infty}\frac{{\pibar_X(y)}}{\Vbar(y)}=\gk^2(0,-\ga )
\whk(0,\ga)
\]
by (\ref{KP}), together with (4.4) and Proposition 5.3 of \cite{kkm},
we thus have
\[
\frac{h(y)}{\Vbar(y)}\to\gk^2(0,-\ga)\whk(0,\ga) \int
_{x\ge
0}\int_{t\ge0} f(x,t)\ga
e^{-\ga x}\,dx\int_{v\ge0} e^{-\ga v} \whV(dt,dv).
\]
Hence, by (\ref{WH}) and Proposition \ref{conconSa},
\begin{eqnarray*}
&& E^{(u)}f \bigl(X_{\tau(u)}-u, \tau(u)- G_{\tau(u)-}\bigr)
\\
&&\qquad \to\frac{ -\Psi_X(-\rmi\ga)}{q} \int_{x\ge0}\int_{t\ge0}
f(x,t)\ga e^{-\ga x}\,dx\int_{v\ge0} e^{-\ga v}
\whV(dt,dv)
\\
&&\hspace*{2pt}\quad\qquad{} +\int_{y\ge0}\frac{\ga}{q}e^{\ga y}\,dy
\int_{{\ceb
}}F(y,\gep) \nbar\bigl(d\gep, \gep(\gz)>y\bigr)+
\rmd_H\frac{\ga}q f(0,0)
\\
&&\qquad = \frac{ -\Psi_X(-\rmi\ga)}{q} \int_{x\ge0}\int_{t\ge0}
f(x,t)\ga e^{-\ga x}\,dx\int_{v\ge0} e^{-\ga v}
\whV(dt,dv)
\\
&&\quad\qquad {} +\frac{\ga}q\int_{x>0}\int
_{t\ge0}f(x,t)\int_{y\ge0} e^{\ga
y}
\Pi_{L^{-1},H}(dt,y+dx) \,dy + \rmd_H\frac{\ga}q f(0,0)
\end{eqnarray*}
by (\ref{Ffxt}). This proves (\ref{Safconw}) since the integral over
$\{
x=0\}$ in the final expression vanishes.
\end{pf}

From (\ref{ltovtS}), a simple calculation shows that the limiting
distribution of the overshoot is as given in (\ref{oSa}), and
%
\begin{eqnarray}
\label{uSa} %
&& P^{(u)}\bigl(\tau(u)-G_ {{\tau(u)}-}\in dt\bigr)
\nonumber\\[-8pt]\\[-8pt]\nonumber
&&\qquad \stackrel{\mathrm {w}} {
\longrightarrow} q^{-1} \bigl(-\Psi_X(-\rmi\ga)\delta
_{-\alpha}^{\whV}(dt)
 +\alpha\rmd_H\delta_0(dt)+K(dt)
\bigr),
\end{eqnarray}
where $K(dt)$ is given by (\ref{K}) and
\[
\delta_{-\alpha}^{\whV}(dt)=\int_{v\ge0}e^{-\alpha v}
\whV(dt,dv).
\]

Using (\ref{Safconw}), we can calculate the limiting value of an EDPF
similar to (\ref{EDPF}); for any $\gb<\ga$ and $\delta\ge0$,
\begin{eqnarray*}
\label{EDPFSa}
&&E^{(u)}e^{\beta(X_{\tau(u)}-u)-\delta(\tau(u)- G_{\tau(u)-})}
\\
&&\qquad\to\frac{ -\Psi_X(-\rmi\ga) }q \int_{x\ge0} \ga
e^{-(\ga
-\gb)
x}\,dx\int_{t\ge0}\int_{v\ge0}
e^{-\gd t-\ga v} \whV(dt,dv)
\\
&&\hspace*{2pt}\quad\qquad{} +\frac{\ga}q\int_{x\ge0}\int
_{t\ge0} e^{\gb x-\gd
t}\int_{y\ge0}
e^{\ga y}\Pi_{L^{-1},H}(dt,y+dx) \,dy + \rmd_H
\frac
{\ga
}q f(0,0)
\\
&&\qquad=\frac{-\ga\Psi_X(-\rmi\ga) }{q(\ga-\gb)\whk(\gd,\ga
)}+\frac
{\ga(\gk(\gd,-\beta)-\gk(\gd,-\ga))}{q(\ga-\beta)}
\end{eqnarray*}
by the same calculation as (\ref{EDPF}).

The results of this section, in the convolution equivalent case, can be
derived from a
path decomposition for the limiting process given in \cite{GM2}. The
main result in \cite{GM2}, Theorem 3.1, makes precise the idea that
under $P^{(u)}$ for large $u$, $X$ behaves like an Esscher transform of $X$
up to an independent exponential time $\tau$. At this time, the process
makes a large jump into a neighborhood of $u$, and if
$W_t=X_{\tau+t}-u$ then
%
\begin{eqnarray}\label{defW}
P(W\in d w) &=& \gk(0,-\ga) \int_{z\in\R}\alpha
e^{-\alpha z} \Vbar(-z)\,d z P_z\bigl(X\in d w\mid \tau(0) <
\infty\bigr),\nonumber
\\[-4pt]
\eqntext{w\in D,}
\end{eqnarray}
where we set $\Vbar(y)=q^{-1}$ for $y<0$.
Thus, $W$ has the law of $X$ conditioned on $\tau(0)< \infty$
and started with initial distribution
\[
\label{W0} P(W_0\in d z)=\gk(0,-\ga) \alpha e^{-\alpha z}
\Vbar(-z) \,d z,\qquad z\in \R.
\]
In the Cram\'er--Lundberg case, there is no comparable decomposition
for the entire path since there is no ``large jump'' at which to do the
decomposition. One of the aims of this paper is to offer an alternative
approach by describing the path from the time of the last maximum prior
to first passage until the time of first passage. This allows the
limiting distribution of many variables associated with ruin to be
readily calculated.


\setcounter{lemma}{0}
\setcounter{teo}{0}
\setcounter{cor}{0}
\setcounter{prop}{0}
\setcounter{remark}{0}

\setcounter{equation}{0}

\begin{appendix}\label{append}
\section*{Appendix: Completion of the proof of Proposition \texorpdfstring{\protect\ref{2222}}{4.2}
when~$X$ is compound Poisson}\label{s9}

For $\vep>0$, let
\[
X^\vep_t=X_t-\vep t.
\]
If $X$ is compound Poisson, then Proposition \ref{2222} holds for
$X^\vep
$. The aim is then to take limits as $\vep\to0$ and check that (\ref
{ch2}) continues to hold in the limit. We begin with an alternative
characterization of the constants in (\ref{n1d}). Recall the notation
of~(\ref{Xpath}).

\begin{lemma}\label{lema1} Assume $0$ is irregular for $(0,\infty)$, then
\[
\mathrm{d}_{\widehat L^{-1}}n(d\gep)=P(X_{[0,\tau(0)]}\in d\gep).
\]
\end{lemma}

\begin{pf}   If $s\notin G$ set $\gep_s=\Delta$ where
$\Delta$ is a cemetery state. Then $\{(t,\gep_{L_{t-}^{-1}})\dvtx\break  t\ge0, \gep_{L_{t-}^{-1}}\neq \Delta\}$
is a Poisson point process with characteristic measure $dt\otimes
n(d\gep)$.
By construction, $n$ is proportional to the law of the first excursion, thus
%
\begin{equation}
\label{fex} n(d\gep)=\llvert n\rrvert P(X_{[0,\tau(0)]}\in d\gep).
\end{equation}
Now let $\gs=\inf\{t\dvtx \gep_{L_{t-}^{-1}}\neq\Delta\}$. Then $\gs$ is
exponentially distributed with parameter~$\llvert  n\rrvert  $. On the other hand $\gs$
is the time of the first jump of $L^{-1}$ and hence is exponential with
parameter $ p$ given by (\ref{nLT2}). A short calculation using duality
(see,\vadjust{\goodbreak} e.g., the paragraph following (2.7) in \cite{Ch}) shows that if
$0$ is irregular for $(0,\infty)$, then
%
\begin{equation}
\label{pd} p\mathrm{d}_{\widehat L^{-1}}=1.
\end{equation}
Hence, $\llvert  n\rrvert  ^{-1}=\mathrm{d}_{\widehat L^{-1}}$ and the result follows from (\ref{fex}).
\end{pf}

Let $n^\vep$ denote the excursion measure of $X^\vep$, with similar
notation for all other quantities related to $X^\vep$ or $\widehat
{X^\vep
}$. To ease the notational complexity, we will write
\[
\widehat{d}_\vep= \mathrm{d}_{(\widehat{L}^\vep)^{-1}} \quad\mbox {and}\quad\widehat{d}=
\mathrm{d}_{\widehat L^{-1}}.
\]

\begin{lemma}\label{lema11} Assume $0$ is irregular for $(0,\infty)$, then
$\widehat{d}_\vep$ is nondecreasing, and for any $\gd\ge0$
\[
\label{dlim} \widehat{d}_\vep\downarrow\widehat{d}_\gd
\qquad\mbox{as } \vep\downarrow \gd.
\]
\end{lemma}

\begin{pf}  Clearly, for $0\le\gd< \vep$, we have $\tau^\gd(0)\le\tau^\vep(0)$
and $ \tau^\vep(0)\downarrow\tau^\gd(0)$ as $\vep\downarrow\gd
$. Thus,
\[
E\bigl(e^{-\tau^\vep(0)};\tau^\vep(0)<\infty\bigr)\uparrow E
\bigl(e^{-\tau^\gd
(0)};\tau ^\gd(0)<\infty\bigr),
\]
and so from (\ref{nLT2}), $p^\vep\uparrow p^\gd$. Hence, by (\ref{pd}),
$
\widehat{d}_\vep\downarrow\widehat{d}_\gd$.
\end{pf}

\begin{prop}\label{pa2} Assume $X$ is compound Poisson and
$f\dvtx [0,\infty
)^2\to[0,\infty)$ is continuous with compact support. Then
\begin{eqnarray*}
\label{lema21} &&
\int_{t\ge0}\int_{z\ge0}f(t,z)
n^\vep\bigl(\gep(t)\in-dz, \zeta >t\bigr)\,dt
\\
&&\qquad \to
\int_{t\ge0}\int_{z\ge0} f(t,z) n
\bigl(\gep(t)\in-dz, \zeta >t\bigr)\,dt \qquad\mbox{as }\vep\to0.
\end{eqnarray*}
\end{prop}

\begin{pf}   Assume $f$ vanishes for
$t\ge r$. Then
%
\begin{equation}
\label{domf} f\bigl(t,-X_t^\vep\bigr)I\bigl(
\tau^\vep(0)>t\bigr)\le\llVert f\rrVert _\infty I(t\le r).
\end{equation}
Thus, using Lemma \ref{lema1},
\begin{eqnarray*}
&& \int_{t\ge0}\int_{z\ge0} f(t,z)
n^\vep \bigl(\gep(t)\in-dz, \zeta >t\bigr)\,dt
\\
&&\qquad =\widehat d^{-1}_\vep\int_{t\ge0}\int
_{z\ge0} f(t,z) P\bigl(X_t^\vep \in-dz,
\tau^\vep(0)>t\bigr)\,dt
\\
&&\qquad =\widehat d^{-1}_\vep\int_{t=0}^\infty
E\bigl(f\bigl(t,-X_t^\vep\bigr);\tau^\vep
(0)>t\bigr)\,dt
\\
&&\qquad \to\widehat d^{-1} \int_{t=0}^\infty E
\bigl(f(t,-X_t);\tau(0)>t\bigr)\,dt
\\
&&\qquad =\int_{t\ge0}\int_{z\ge0} f(t,z) n\bigl(
\gep(t)\in-dz, \zeta>t\bigr)\,dt
\end{eqnarray*}
by (\ref{domf}) and dominated convergence, since $X^\vep_t\to X_t$,
$\tau^\vep(0)\to\tau(0)$ and $P(\tau(0)=t)=0$.
\end{pf}

\begin{prop}\label{pa3} Assume $X$ is compound Poisson and
$f\dvtx [0,\infty
)^2\to[0,\infty)$ is continuous with compact support. Then
\[
\label{lema22} \int_{t\ge0}\int_{z\ge0}
f(t,z) \widehat V^\vep(dt,dz)\to\int_{t\ge0}\int
_{z\ge0} f(t,z) \widehat V(dt,dz) \qquad\mbox{as }\vep\to0.
\]
\end{prop}

\begin{pf}   We will show
%
\begin{equation}
\label{LHc} \bigl(\widehat{L}^\vep\bigr)^{-1}_{s}
\to\widehat{L}^{-1}_{s},\qquad \widehat {H}^\vep_s
\to\widehat {H}_s\qquad\mbox{for all $s\ge0$ as }\vep\to0,
\end{equation}
and that the family
%
\begin{equation}
\label{fd1} f\bigl(\bigl(\widehat{L}^\vep\bigr)^{-1}_{s},
\widehat{H}^\vep_s\bigr),\qquad 0<\vep\le1,
\end{equation}
is dominated by an integrable function with respect to $P\times ds$. Then
\begin{eqnarray*}
\int_{t\ge0}\int_{z\ge0} f(t,z)
\whV^\vep(dt,dz) &=&\int_{t=0}^\infty Ef
\bigl(\bigl(\widehat{L}^\vep\bigr)^{-1}_{s},
\widehat {H}^\vep_s\bigr)\,ds
\\
&\to&\int_{t=0}^\infty Ef\bigl(\widehat{L}^{-1}_{s},
\widehat{H}_s\bigr)\,ds
\\
&=& \int_{t\ge0}\int_{z\ge0} f(t,z) V(dt,dz).
\end{eqnarray*}

For $\vep\ge0$, let
$
A_\vep=\{s\dvtx  \overline{\widehat{X}}{}^\vep_s=\widehat{X}^\vep_s\}$.
Then for $0\le\gd< \vep$,
$
A_\gd\subset A_\vep$.
Further, for any~$T$, if $\vep$ is sufficiently close to $0$, then
$
A_0\cap[0,T]=A_\vep\cap[0,T]$.
Thus, by Theorem 6.8 and Corollary 6.11 of \cite{kypbook},
\[
\widehat d_\gd\widehat{L}^\gd_t =\int
_0^t I_{A_\gd}(s)\,ds \le\int
_0^t I_{A_\vep}(s)\,ds =\widehat
d_\vep\widehat{L}^\vep_t, \qquad\mbox{all } 0
\le\gd < \vep,
\]
and
\[
\widehat d \widehat{L}_t=\widehat d_\vep
\widehat{L}^\vep_t,\qquad 0\le t\le T,
\]
if $\vep$ is sufficiently close to $0$.
Hence, for all $0\le\gd< \vep$,
%
\begin{equation}
\label{Lgdv} \bigl(\widehat{L}^\gd\bigr)^{-1}_s
=\inf\bigl\{t\dvtx \widehat{L}^\gd_t>s\bigr\} \ge\inf\bigl
\{t\dvtx (\widehat d_\vep/\widehat d_\gd)\widehat{L}^\vep_t>s
\bigr\} =\bigl(\widehat{L}^\vep\bigr)^{-1}_{(\widehat d_\gd/\widehat d_\vep)s},
\end{equation}
with equality if $\gd=0$ and $\vep$ is sufficiently close to $0$.

Fix $s\ge0$ and assume
$\vep$ is sufficiently close to $0$ that equality holds in (\ref
{Lgdv}) with $\gd=0$. Thus,
%
\begin{equation}
\label{hat1} \bigl(\widehat{L}^\vep\bigr)^{-1}_{s}=
\widehat{L}^{-1}_{(\widehat d_\vep
/\widehat d)s}.
\end{equation}
Since
$(\overline{\widehat{X}}{}^\vep)_t
=\overline{\widehat{X}}_t+J_{\vep,t}$
where
%
\begin{equation}
\label{hat2} 0\le J_{\vep,t} \le\vep t,
\end{equation}
it then follows that
%
\begin{equation}
\label{hat3} \widehat{H}^\vep_s=\bigl(\overline{
\widehat{X^\vep}}\bigr)_{(\widehat
{L}^\vep)^{-1}_{s}} =\overline{
\widehat{X}}_{\widehat{L}^{-1}_{(\widehat d_\vep
/\widehat d)s}}+J_{\vep, \widehat
{L}^{-1}_{(\widehat d_\vep/\widehat d)s}} =\widehat{H}_{(\widehat d_\vep/\widehat d)s}+J_{\vep, \widehat
{L}^{-1}_{(\widehat d_\vep/\widehat d)s}}.
\end{equation}
Hence, using Lemma \ref{lema11}, (\ref{LHc}) follows from (\ref{hat1}),
(\ref{hat2}) and (\ref{hat3}).

Now let $0\le\vep\le1$. Then by (\ref{Lgdv}),
\[
\bigl(\widehat{L}^\vep\bigr)^{-1}_{s}\ge\bigl(
\widehat{L}^1\bigr)^{-1}_{(\widehat
d_\vep/\widehat d_1)s}.
\]
Thus, by monotonicity of $\widehat d_\vep$,
\[
I\bigl(\bigl(\widehat{L}^\vep\bigr)^{-1}_{s}\le
r\bigr) \le I\bigl(\bigl(\widehat{L}^1\bigr)^{-1}_{(\widehat d_\vep/\widehat d_1)s}
\le r\bigr) \le I\bigl(\bigl(\widehat{L}^1\bigr)^{-1}_{(\widehat d/\widehat d_1)s}
\le r\bigr).
\]
Hence, if $f$ vanishes for
$t\ge r$, then
\[
f\bigl(\bigl(\widehat{L}^\vep\bigr)^{-1}_{s},
\widehat{H}^\vep_s\bigr)\le\llVert f\rrVert
_\infty I\bigl(\bigl(\widehat {L}^1\bigr)^{-1}_{(\widehat d/\widehat d_1)s}
\le r\bigr),
\]
where
\[
E\int_0^\infty I\bigl(\bigl(
\widehat{L}^1\bigr)^{-1}_{(\widehat d/\widehat d_1)s}\le r\bigr) \,ds = (
\widehat d_1/\widehat d) E\int_0^\infty
I\bigl(\bigl(\widehat {L}^1\bigr)^{-1}_{s}\le
r\bigr) \,ds<\infty,
\]
which proves (\ref{fd1}).
\end{pf}

\begin{pf*}{Proof  of Proposition \ref{2222} when $X$ is compound Poisson}
Assume $X$ is compound Poisson. Since $\rd_{L^{-1}}=0$ whenever $0$ is
irregular for $(0,\infty)$, it follows that $\mathrm{d}_{(L^\ve
)^{-1}}=\rd_{L^{-1}}=0$. Further, (\ref{ch2}) holds for $X^{\ve}$. Hence,
(\ref{ch2}) for $X$ follows from Propositions \ref{pa2} and \ref{pa3}.
\end{pf*}
\end{appendix}

\section*{Acknowledgement}
I would like to thank Ron Doney for his
help with parts of Section~\ref{s3}, and the referees for their useful
suggestions.



\printaddresses
\end{document}